\documentclass[12pt, a4paper]{article}
\usepackage[left=1in,right=1in,top=1.0in,bottom=1in]{geometry}

\usepackage{amsmath,amssymb,amsthm,amsfonts,afterpage,hyperref}
\usepackage{amscd,subeqnarray}

\usepackage{graphicx}
\usepackage{caption}
\usepackage{subcaption}
\usepackage{authblk}
\usepackage{tikz}

\usepackage{array}
\usepackage{wrapfig}
\usepackage{color}
\usepackage{ragged2e}
\usepackage{stmaryrd}
\usepackage{siunitx}
\usepackage{multirow}
\usepackage{xcolor,cancel}
\usepackage{boldfonts}
\usepackage{symbols}
\usepackage{footmisc}
\usepackage{float}
\usepackage{setspace}
\usepackage{bm}
\usepackage{booktabs}
\usepackage{tikz, xcolor}
\usepackage{soul}

\usepackage[T1]{fontenc}
\usepackage[utf8]{inputenc}
\usepackage{mathrsfs}

\usepackage[linesnumbered,ruled,vlined]{algorithm2e}
\SetKwInput{KwInput}{Input}                
\SetKwInput{KwOutput}{Output} 

\usepackage[sort, numbers]{natbib}
\setcitestyle{square}

\usetikzlibrary{shapes,arrows, positioning,calc}	
\usetikzlibrary{arrows,snakes,backgrounds}

\newcommand{\bfa}[1]{\boldsymbol{#1}} 			
			
\newcommand{\bfeps}{\boldsymbol{\epsilon}}

\newcommand{\bfsigma}{\boldsymbol{\sigma}}

\newcommand{\Sym}{\text{Sym}}   			%
\newcommand{\tr}{\text{tr}}       				%

\DeclareMathAlphabet{\mathpzc}{OT1}{pzc}{m}{it}

\newcommand{\bfu}{\boldsymbol{u}}

\newcommand{\bfx}{\boldsymbol{x}}

\newcommand{\bfT}{\boldsymbol{T}}		
\newcommand{\bfzero}{\boldsymbol{0}}

\newcommand{\bff}{\boldsymbol{f}}	
	
\newtheorem{theorem}{Theorem}[section]

\newtheorem{remark}{Remark}

\newtheorem*{cwf}{Continuous weak formulation}
\newtheorem*{dwf}{Discrete weak formulation}

\providecommand{\keywords}[1]
{
  \small	
  \textbf{\textit{Keywords---}} #1
}

\title{Computational investigation of crack-tip fields in a compressed nonlinear strain-limiting material\footnote{This work is dedicated to the memory of Professor K. R. Rajagopal (Former Distinguished Professor of Mechanical Engineering, Texas A\&M University, College Station, Texas, USA), whose pioneering contributions to the field of solid mechanics were foundational to this research and a deep influence on the authors.}}

\author[1]{Dambaru Bhatta}
\author[2]{Saugata Ghosh}
\author[3]{S. M. Mallikarjunaiah} 

\affil[1,2]{School of Mathematical \& Statistical Sciences,\\
The University of Texas - Rio Grande Valley, \\
Edinburg, Texas 78539, USA \\
Email:  dambaru.bhatta@utrgv.edu, saugata.ghosh01@utrgv.edu}

\affil[3]{Department of Mathematics \& Statistics,\\
Texas A\&M University-Corpus Christi, \\
Corpus Christi, Texas 78412-5825, USA \\
Email: m.muddamallappa@tamucc.edu (Corresponding Author)}
\date{}

\begin{document}

\maketitle
	    
\begin{abstract}
A finite element framework is presented for the analysis of crack-tip phenomena in an elastic material containing a single edge crack under compressive loading. The mechanical response of the material is modeled by a nonlinear constitutive relationship that algebraically relates stress to linearized strain. This approach serves to mitigate non-physical strain singularities and ensures that the crack-tip strains don't grow, unlike singular stresses.  A significant advancement is thus achieved in the formulation of boundary value problems (BVPs) for such complex scenarios. The governing equilibrium equation, derived from the balance of linear momentum and the nonlinear constitutive model, is formulated as a second-order, vector-valued, quasilinear elliptic BVP. A classical traction-free boundary condition is imposed on the crack face. The problem is solved using a robust numerical scheme in which a Picard-type linearization is combined with a continuous Galerkin finite element method for the discretization. Analyses are performed for both an isotropic and a transversely isotropic elastic solid containing a crack subjected to compressive loads. The primary crack-tip variables—stress, strain, and strain energy density—are examined in detail. It is demonstrated that while high concentrations of compressive stress and strain energy density are observed at the crack tip, the growth of strain is substantially lower than that of stress. These findings are shown to be consistent with the predictions of linear elastic fracture mechanics, but a more physically meaningful representation of the crack-tip field is provided by the nonlinear approach. A rigorous basis is thus established for investigating fundamental processes like crack propagation and damage in anisotropic, strain-limiting solids under various loading conditions, including compression.
 \end{abstract}
	 
\noindent \keywords{Finite element method; Compressive load; Transversely Isotropic; Nonlinear constitutive relations}

\section{Introduction}\label{Introduction}
A comprehensive understanding of crack behavior under compressive loads is of significant engineering and medical importance, transcending traditional fracture mechanics focused on tensile failure. In engineering, this knowledge is critical for the safety and reliability of large-scale structures built from quasi-brittle materials like concrete and rock. Civil infrastructure, such as bridges, tunnels, and dams, is constantly subjected to immense compressive forces, and the propagation of existing cracks in these materials can lead to catastrophic failure modes like buckling and progressive collapse \cite{kim2025effect,bahmani2024effect}. Furthermore, in the aerospace and marine sectors, components such as stiffened panels in ships and aircraft fuselages experience complex cyclic compressive loads, where crack-face contact and friction can dramatically alter the stress distribution and structural integrity \cite{khedmati2009sensitivity}. The development of accurate numerical models is essential for predicting material failure and assessing the residual strength of these vital structures.

In the medical field, this research is equally indispensable for improving patient outcomes. The study of compressive fracture mechanics is fundamental to understanding and preventing vertebral compression fractures, a common and debilitating condition in aging populations with osteoporosis \cite{parreira2017overview}. By analyzing microcrack behavior in bone tissue under compression, we can gain insights into the mechanisms of bone remodeling and repair, which is crucial for developing new treatments and therapies \cite{shu2014surface}. This field of study also directly informs the design of more durable and biocompatible orthopedic implants and dental prosthetics, as these devices must endure significant cyclic compressive forces from bodily movements and mastication \cite{manfredini2024clinical}. Ultimately, a deeper understanding of compressive crack phenomena provides a rigorous foundation for developing advanced materials and therapeutic strategies that enhance structural longevity and human health.

A comprehensive understanding of crack-tip phenomena under compressive loading is crucial for a variety of engineering and scientific disciplines. Unlike simple tensile loading, which is well-characterized by classical linear elastic fracture mechanics, compressive forces introduce complex behaviors such as crack-face contact and friction. These interactions significantly alter the local stress, strain, and energy fields around the crack tip, necessitating a detailed investigation into these parameters for both isotropic and transversely isotropic materials. While studying crack-tip stresses is essential for identifying potential failure initiation sites, analyzing crack-tip strains provides critical information about the material's local deformation and the onset of damage \cite{broberg1999}. The strain energy density criterion, which synthesizes both stress and strain, offers a more robust framework for predicting crack initiation and propagation under such intricate conditions \cite{sih1973handbook, yoon2021quasi}. Historically, the analytical foundation for characterizing stress concentrations has been rooted in linearized elasticity theory, which, despite its utility, predicts physically unrealistic strain singularities at the crack tip \cite{Inglis1913, love2013treatise}. This discrepancy has motivated extensive efforts to develop more physically representative constitutive models that can effectively mitigate these singularities and provide a more accurate representation of material behavior, particularly under compressive loads where a physically meaningful theory is essential for structural integrity. Despite its widespread use for modeling crack initiation and propagation, linear elastic fracture mechanics (LEFM) has inherent limitations that extend beyond the well-known strain singularity. A major issue is that LEFM predicts a physically unrealistic blunt crack-opening profile and the interpenetration of crack faces, especially in interfaces between different materials. The problem of crack-tip singularity isn't fully resolved even in various nonlinear elasticity models. For instance, the work of Knowles and Sternberg \cite{knowles1983} and models incorporating specific constraints, such as the Bell constraint model \cite{tarantino1997}, have highlighted that these nonlinear frameworks don't completely eliminate the singularity. These persistent issues raise a fundamental question that drives the present research: can algebraic nonlinear constitutive models effectively mitigate the crack-tip strain singularity, even if some stress singularity remains?

Classical elasticity theories, such as those formulated by Cauchy and Green, often fall short when describing material behavior under extreme deformations, particularly near crack tips, where they predict physically impossible singularities. In response to these limitations, Rajagopal and his collaborators have developed a generalized framework for elasticity \cite{rajagopal2003implicit,rajagopal2007elasticity,rajagopal2007response,rajagopal2009class,rajagopal2011non,rajagopal2011conspectus,rajagopal2014nonlinear,rajagopal2016novel,rajagopal2018note} that is rooted in robust thermodynamic principles. This novel approach, which uses \textit{implicit constitutive models}, defines an elastic body as non-dissipative and characterized by a relationship that implicitly links the Cauchy stress and deformation gradient tensors \cite{bustamante2018nonlinear,bustamante2021new,bustamante2015implicit,rodriguez2021stretch}. A key advantage of this framework is its ability to produce a hierarchy of explicit nonlinear relationships, allowing strain to be expressed as a nonlinear function of stress. A specific subset of these models is particularly powerful because it exhibits a unique \textit{strain-limiting property}.  This means the models can represent linearized strain as a uniformly bounded function throughout the material, even when subjected to high levels of stress. This property makes these models exceptionally well-suited for investigating crack and fracture behavior in brittle materials \cite{rajagopal2011modeling, gou2015modeling, Mallikarjunaiah2015, MalliPhD2015}. The utility of these strain-limiting models has been demonstrated through numerous studies that have revisited and provided new insights into classical elasticity problems \cite{kulvait2013,rajagopal2018bodies,bustamante2009some,bulivcek2014elastic, erbay2015traveling, zhu2016nonlinear, csengul2018viscoelasticity, itou2018states, itou2017contacting, yoon2022CNSNS, yoon2022MMS,rodriguez2022longitudinal,bustamante2025piecewise,rajagopal2024mathematical}. Their versatility in explaining the mechanical behavior of a wide range of materials, especially concerning crack and fracture phenomena, is a significant advantage. For example, recent research has shown that when formulated within this strain-limiting framework, quasi-static crack evolution problems can predict complex crack patterns and even increased crack-tip propagation velocities \cite{lee2022finite, yoon2021quasi}.

This study presents a rigorous analysis of a singular crack within an isotropic and transversely isotropic solid under compressive load. The theoretical foundation is built upon a novel, algebraically nonlinear constitutive model specifically formulated to describe the stress-strain behavior of isotropic and orthotropic materials under compressive forces. By coupling this constitutive law with the principle of linear momentum balance, we formulate the problem as a vector-valued, quasi-linear elliptic boundary value problem. This mathematical framework, which avoids the unphysical singularities of classical linear elasticity, serves as the core idea to revisit the classical problems in solid mechanics and fracture mechanics.  Given the analytical intractability of the resulting system of nonlinear partial differential equations, we developed a computational solution strategy. The domain is discretized using the finite element method, a technique chosen for its proven capability to accurately capture field behavior in the vicinity of geometric singularities. To address the inherent nonlinearities of the governing equations, we implemented a robust Picard iterative algorithm. Numerical convergence of the solution is rigorously established by monitoring the progressive reduction of the residual norm throughout the iterative process. Our computational results elucidate several key phenomena. The simulations confirm a controlled and bounded growth of the crack-tip strains, effectively regularizing the singularity predicted by classical theory. A parametric study reveals that the material's response to compressive loading is critically dependent on one of the two constitutive parameters. For a specific range of this parameter, we observe significant stress concentration and a notable decrease in strain-energy density around the crack tip. Conversely, a distinct, contrary behavior is observed for other parameter values. These foundational findings open several promising avenues for future research, including the extension of this framework to quasi-static thermo-elastic cracks and the intricate dynamics of crack propagation in transversely isotropic materials. A particularly compelling direction involves a dedicated investigation into the behavior of these cracks under purely compressive loads, which pose unique mathematical challenges and are known to induce distinct failure modes, such as buckling and crushing, not observed in tension.

The structure of the paper is as follows. 
Section~\ref{math_formulation} establishes the theoretical groundwork by introducing the foundational implicit theory and deriving the nonlinear constitutive relation. 
Following this, Section~\ref{sec:BVP} formulates the specific boundary value problem for a static crack and demonstrates the existence of a unique solution to its weak form. 
Our numerical scheme, which combines a continuous Galerkin finite element method with Picard's iterative algorithm, is detailed in Section~\ref{sec:fem}. 
Section~\ref{sec:rd} is dedicated to presenting and analyzing the numerical results, with a particular focus on the effects of key model parameters on the compressive crack-tip field variables. 
The paper concludes by summarizing the principal findings and their significance, along with some brief outlines of future work within the context of the model studied in this paper.

\section{Mathematical Preliminaries} \label{math_formulation}

Let $\Lambda \subset \mathbb{R}^2$ be a closed, bounded, and connected Lipschitz domain representing the reference configuration of a material body. The boundary of this domain, $\partial \Lambda$, is partitioned into two disjoint parts: a non-empty Dirichlet boundary $\Gamma_D$ and a Neumann boundary $\Gamma_N$, such that $\partial \Lambda = \overline{\Gamma_D} \cup \overline{\Gamma_N}$ and $\Gamma_D \cap \Gamma_N = \emptyset$. We denote a point in the reference configuration by $\bm{X} = (X, \, Y)$ and its corresponding position in the deformed configuration by $\bm{x} = (x, \, y)$. The deformation is described by the displacement field $\bm{u}: \Lambda \to \mathbb{R}^2$, where 
\[
\bm{u}(\bm{X}) = \bm{x}(\bm{X}) - \bm{X}.
\]

For the subsequent analysis, we adopt the standard notations of Lebesgue and Sobolev spaces \cite{ciarlet2002finite, evans1998partial}. The space $L^p(\Lambda)$ for $p \in [1, \infty)$ denotes the set of all Lebesgue integrable functions on $\Lambda$. Specifically, $L^2(\Lambda)$ is the space of square-integrable functions, equipped with the inner product $\langle v, w \rangle := \int_{\Lambda} v w \, d\bm{x}$ and the induced norm $\|v\|_{L^2(\Lambda)} := \sqrt{\langle v, v \rangle}$. The Sobolev space $H^1(\Lambda) := W^{1,2}(\Lambda)$ consists of functions in $L^2(\Lambda)$ whose weak derivatives are also in $L^2(\Lambda)$, and it is equipped with its standard inner product and norm. We define $H^1_0(\Lambda)$ as the closure of the space of compactly supported smooth functions, $C_0^\infty(\Lambda)$, under the $H^1$-norm.

We define the space of symmetric second-order tensors in $\mathbb{R}^{2 \times 2}$ as $\mathrm{Sym}(\mathbb{R}^{2 \times 2})$. This space is endowed with the inner product $\langle \bm{A}, \bm{B} \rangle = \sum_{i,j=1}^2 A_{ij} B_{ij}$ and the associated norm $\|\bm{A}\| = \sqrt{\langle \bm{A}, \bm{A} \rangle}$. The kinematics of the deformation are described by the deformation gradient 
\[
\bm{F} = \nabla_{\bm{X}} \bm{x} = \bm{I} + \nabla \bm{u}.
\]
Under the assumption of small strains, where $\|\nabla\bm{u}\|$ is of a lower order of magnitude, the standard kinematic tensors simplify to their linearized forms. Specifically, the linearized strain tensor is given by 
\[
\bm{\epsilon} = \frac{1}{2}(\nabla\bm{u} + \nabla\bm{u}^T),
\]
 and its relationship to the full nonlinear kinematic measures becomes approximately linear. The Cauchy stress tensor $\bfsigma$ satisfies the balance of linear momentum, and it is symmetric. 

\subsection{Rajagopal's theory of elasticity and strain-limiting models}

Unlike classical Cauchy elasticity, which assumes an explicit relationship 
\begin{equation}
\bfsigma =\widehat{\mathcal{F}}(\bm{F}), 
\end{equation}
our approach is founded on the generalized theory of implicit elasticity pioneered by Rajagopal  \cite{rajagopal2003implicit,rajagopal2007elasticity,rajagopal2011non,rajagopal2011conspectus,rajagopal2014nonlinear,rajagopal2007response}. In this framework, the constitutive response is defined implicitly. A particular subclass of these models is the \textit{strain-limiting constitutive theories}, which specify the linearized strain as a uniformly bounded function of the stress tensor, and such a notion, within the first-order assumption,  is expressed in the form
\begin{equation}
 \bm{\epsilon} = \mathcal{F}(\bfsigma), 
\end{equation} 
 where $\mathcal{F}$ is a tensor-valued function stress tensor $\bfsigma$ and its invariants. The key feature of such a model is that, 
 \[
 \|\bm{\epsilon}(\bm{u})\| \leq M,
 \]
 for some constant $M > 0$, regardless of the magnitude of the stress. This strain-limiting property is crucial for mitigating the unphysical strain singularities predicted at crack tips by classical linear elastic fracture mechanics. The specific constitutive relationship used in this work is an algebraically nonlinear form of this model, given by:
\begin{equation}
\bm{\epsilon} = \Psi_0(\mathrm{tr} \, \bfsigma, \|\bfsigma\|) \bm{I} + \Psi_1(\|\bfsigma\|) \bfsigma
\end{equation}
where $\Psi_0$ and $\Psi_1$ are scalar functions of the principal invariants of the stress tensor. This constitutive law is combined with the equilibrium equations and compatibility conditions:
\[
 \mathrm{curl} \mathrm{curl} \, \bm{\epsilon} = \bm{0},
 \]
 to form a complete system of equations governing the problem. This novel approach allows us to investigate issues involving materials containing cracks and fractures while maintaining physically realistic, bounded strain fields throughout the domain.

In this work, we adopt a specific subclass of \textit{Rajagopal’s theory of elasticity} known as \textit{strain-limiting constitutive models}. A key mathematical property of these models is that the strain response is uniformly bounded. This means that for a given tensor-valued function $\mathcal{F}(\cdot)$ relating stress to strain, there exists a constant $M>0$ such that $\|\mathcal{F}(\bfsigma)\| \leq M$ for all symmetric stress tensors $\bfsigma \in \mathrm{Sym}(\mathbb{R}^{2 \times 2})$. This uniform bound is crucial for circumventing the unphysical strain singularities predicted by classical linear elastic fracture mechanics at crack tips.

For our numerical study, we define the specific functional form of the strain-limiting relationship. The linearized strain tensor $\bm{\epsilon}$ is related to the Cauchy stress tensor $\bfsigma$ via the constitutive law:
\begin{equation}
\bm{\epsilon} = \mathcal{F}(\bfsigma) :=  \Psi(\|\bfsigma\|) \, \bfsigma
\end{equation}
where  $\Psi: \mathbb{R}_+ \to \mathbb{R}$ is a scalar function of the stress invariants. Following prior work in this field \cite{yoon2021quasi,yoon2022MMS,yoon2022CNSNS,lee2022finite,itou2017contacting,itou2017nonlinear,itou2018states}, we employ a particular form of the function $\mathcal{F}(\cdot)$ that guarantees a strain-limiting response:
\begin{equation}\label{def-F}
\mathcal{F}(\bfsigma) = \frac{\mathbb{K}[\bfsigma]}{\left( 1 + \beta^{\alpha} \|\mathbb{K}^{1/2}[\bfsigma]\|^{\alpha} \right)^{1/\alpha}}, \quad \beta \geq 0, \alpha > 0
\end{equation}
This form ensures that the strain response is uniformly bounded, with a maximum value of $1/\beta$, i.e., 
\[
\sup_{\bfsigma \in \mathrm{Sym}} \|\mathcal{F}(\bfsigma)\| \leq 1/\beta.
\]
 It is also important to note that as the model parameters approach certain limits---specifically as $\beta \to 0$ or $\alpha \to \infty$---this nonlinear strain-limiting model converges to the classical linearized elastic model.  The fourth-order tensor $\mathbb{E}[\cdot]$ represents the \textit{linearized elasticity tensor}. For an orthotropic material, it is defined as:
\begin{subequations}\label{E-def}
\begin{align}
\mathbb{E}[\bm{\epsilon}] &:= 2\mu\bfeps + \lambda\,\tr(\bfeps)\,\bfa{I} + \gamma(\bfeps \colon \bfa{M})\,\bfa{M}. 
\end{align}
\end{subequations}
Here, $\mu > 0$ and $\lambda > 0$ are the Lam\'e coefficients, $\gamma$ is an additional modulus characterizing the anisotropic response, and the structural tensor $\bfa{M}$ defines the preferred material direction (e.g., fiber orientation) \cite{Mallikarjunaiah2015,ghosh2025finite,ghosh2025computational}. Another fourth-order tensor $\mathbb{K}[\cdot]$ is the inverse of the elasticity tensor, the compliance tensor.  Both $\mathbb{K}[\cdot]$ and $\mathbb{E}[\cdot]$ are positive definite, self-adjoint operators on $\mathrm{Sym}(\mathbb{R}^{2 \times 2})$.

A detailed examination of this tensor-valued function $\mathcal{F}(\bfsigma)$, as described in \cite{Mallikarjunaiah2015,itou2017nonlinear,itou2018states,vasilyeva2024generalized,ghosh2025computational,ghosh2025finite}, reveals several important mathematical properties that underpin the well-posedness of our problem:
\begin{itemize}
    \item \textbf{Uniform Boundedness:} The function $\mathcal{F}$ is uniformly bounded, meaning $\|\mathcal{F}(\bfsigma)\| \leq M$ for all $\bfsigma$. This property guarantees that the strain fields remain finite throughout the domain, even in regions of high stress concentration.
    \item \textbf{Strict Monotonicity:} The function is strictly monotone in the sense that $\langle \mathcal{F}(\bfsigma_1) - \mathcal{F}(\bfsigma_2), \bfsigma_1 - \bfsigma_2 \rangle > 0$ for any distinct tensors $\bfsigma_1, \bfsigma_2 \in \mathrm{Sym}(\mathbb{R}^{2 \times 2})$. This property is crucial for the uniqueness and stability of the solution.
    \item \textbf{Continuity (Lipschitz):} The function is continuous, satisfying a Lipschitz-like condition: $\langle \mathcal{F}(\bfsigma_1) - \mathcal{F}(\bfsigma_2), \bfsigma_1 - \bfsigma_2 \rangle \leq \widehat{c}_1 \|\bfsigma_1 - \bfsigma_2\|^2$. This continuity is essential for proving the existence of a solution.
\item \textbf{Coercivity:} The response function $\mathcal{F}(\cdot)$ is coercive, which is a crucial property for guaranteeing a unique solution to the problem's variational formulation. Specifically, for any symmetric tensor $\bfa{\Pi} \in \Sym(\mathbb{R}^{2 \times 2})$ and any non-zero vector $\bfa{v} \in \mathbb{R}^{2}$, there exists a constant $\widehat{c}_2 > 0$ such that the inequality $ {\left|  \bfa{v} \cdot \mathcal{F}(\bfa{\Pi}) \bfa{v} \right| \geq \widehat{c}_2 \| \bfa{v} \|^2  }$ holds. This constant $\widehat{c}_2$ is a function of the model's parameters, including the Lam\'{e} parameters, and the dimensionality of the problem.
\end{itemize}

The properties of the response function used in this research have been extensively studied. A particularly important point, which has been examined by several researchers, is the existence of hyperelastic models within Rajagopal's theory. This concept is crucial for our work, particularly for formulating the hyperelastic response function that we will use in our numerical simulations.

\begin{remark}
For sufficiently small values of the nonlinearity parameter $\beta$, the constitutive relation is invertible, a fact established in \cite{mai2015strong,mai2015monotonicity}. This inverted relationship can be derived from a scalar potential, meaning the material model is {hyperelastic}. This hyperelastic form, which expresses stress as a function of strain, offers significant computational advantages for displacement-based finite element methods. The specific form of this relationship is given by:
\begin{equation}
\bfsigma(\bfeps) := \Psi\left( \| \mathbb{E}^{1/2}[\bfeps] \| \right) \, \mathbb{E}[\bfeps], \quad \text{where} \quad \Psi(s) = \frac{1}{\left(1 - (\beta s)^{\alpha} \right)^{1/\alpha}}. \label{eq:inverted_hyperelastic}
\end{equation}
In the subsequent analysis, we will use this hyperelastic formulation to construct the governing boundary value problem. A primary goal of our research is to systematically compare the mechanical fields (stress, strain, and strain energy density) at the crack tip predicted by this nonlinear, strain-limiting model against the predictions of its classical linear elastic counterpart, which is recovered by setting $\beta = 0$.
\end{remark}

\begin{remark}
The general class of implicit constitutive theories pioneered by Rajagopal typically does not permit an analytical inversion to express the stress as an explicit function of the strain. As noted in the literature \cite{rajagopal2018note,sfyris2024hyperbolicity}, this non-invertibility is a fundamental feature of these advanced models. In contrast, the specific constitutive relation proposed in this work is deliberately constructed to be invertible. This special characteristic allows us to derive a corresponding explicit hyperelastic model, where the stress can be obtained directly from a strain-energy potential. We wish to emphasize that this inversion is performed for a purely numerical purpose. Its utility lies not in the theoretical framework itself, but in practical applications such as simplifying aspects of the computational implementation.
\end{remark}

\section{Governing equations and well-posedness}
\label{sec:BVP}

The analysis of fracture mechanics in transversely isotropic materials is a cornerstone of modern continuum mechanics, with profound implications for both materials science and structural engineering. The critical importance of this field stems from the widespread use of such materials in high-performance and safety-critical applications. This class of materials, characterized by directionally dependent mechanical properties, includes advanced composites, natural substances like wood, and various geological formations and biological tissues. A profound understanding of the initiation and propagation of cracks is essential, as these can severely compromise a structure's integrity and lead to catastrophic failure. For example, in the aerospace industry, where fiber-reinforced composites are ubiquitous, the ability to accurately predict and manage crack growth is fundamental to ensuring the safety and reliability of aircraft. Similarly, in civil engineering, a precise comprehension of fracture processes is vital for assessing the long-term durability of infrastructure.

This section formally establishes the governing boundary value problem (BVP) for a cracked, transversely isotropic solid that exhibits strain-limiting behavior. We begin by outlining the key assumptions necessary to ensure the mathematical well-posedness of the problem and then present the weak (variational) formulation. This weak formulation provides a rigorous foundation for both the theoretical existence proof of a solution and the subsequent numerical approximation using methods such as the finite element method. The mechanical system is governed by the principle of conservation of linear momentum, which, in the absence of inertial effects, simplifies to the static equilibrium equation. This principle states that the divergence of the Cauchy stress tensor $\bfsigma$ is balanced by the body force vector $\bff$. This equilibrium equation is fundamentally coupled with the hyperelastic constitutive relationship, which expresses the stress as a nonlinear function of the strain, i.e., $\bfsigma = \bfsigma(\bfeps)$. The system is closed by imposing essential (Dirichlet) and natural (Neumann) boundary conditions on complementary parts of the domain's boundary, $\partial \Lambda$. The resulting strong form of the BVP is to find the displacement field $\bfa{u}$ that satisfies:
\begin{subequations}
\label{eq:strong_form_bvp}
\begin{align}
-\nabla \cdot \bfsigma &= \bfa{f},  && \text{in } \Lambda, \label{eq:equilibrium} \\
\bfa{u} &= \bfa{u}_0, && \text{on } \Gamma_D, \label{eq:dirichlet_bc} \\
\bfsigma \bfa{n} &= \bfa{g}, && \text{on } \Gamma_N, \label{eq:neumann_bc}
\end{align}
\end{subequations}
where $\bfa{f} \in (L^2(\Lambda))^2$ represents the body force per unit volume and $\bfa{g} \in (L^{2}(\Gamma_N))^2$ is the prescribed traction vector on the Neumann boundary $\Gamma_N$.

For the BVP described by \eqref{eq:strong_form_bvp} to be mathematically well-posed, ensuring a unique and stable solution, we introduce the following physically motivated assumptions:

\begin{itemize}
 \item[\textbf{A1:}] \textbf{Model and material parameters.} The scalar modeling parameters, $\alpha$ and $\beta$, which control the nonlinear constitutive response, are assumed to be positive and uniform constants throughout the domain. Similarly, the Lamé parameters $\mu$ and $\lambda$ are also taken as positive, uniform constants, consistent with a homogeneous material. A crucial validation of our model involves numerically verifying that in the limit as $\beta \to 0^+$, the predictions of the nonlinear model converge to those of classical linear elasticity.
  \item[\textbf{A2:}] \textbf{Static equilibrium condition.} For the pure traction problem where the displacement is unconstrained ($\Gamma_D = \emptyset$), the external loads must be self-equilibrated. This physical necessity imposes the following integral compatibility condition on the body forces and surface tractions:
   \begin{equation}
    \int_{\Lambda} \bfa{f} \, d\bfa{x} + \int_{\partial \Lambda} \bfa{g} \, ds = \bfa{0}.
    \end{equation}
 This condition ensures that the net force acting on the body is zero, thereby precluding any rigid body motion.
\item[\textbf{A3:}] \textbf{Regularity of boundary data.} The prescribed displacement on the Dirichlet boundary, $\bfa{u}_0$, is assumed to possess sufficient regularity, typically $\bfa{u}_0 \in (W^{1,1}(\Lambda))^2$. This ensures that the boundary deformation is physically reasonable and the problem remains mathematically tractable, as this regularity is a prerequisite for a well-defined trace on the boundary.
\end{itemize}

To construct the weak formulation of the boundary value problem and establish the existence of a solution, we first define two crucial subspaces of the Sobolev space $\left( H^1(\Lambda)\right)^2$. These spaces are essential for distinguishing between the kinematically admissible variations (test functions) and the space of all possible solutions (trial functions).

The space for test functions, $\bfa{V}_{\bfzero}$, consists of all displacement fields from the Sobolev space $\left( H^1(\Lambda)\right)^2$ that vanish on the Dirichlet boundary $\Gamma_D$. This condition ensures that the test functions are compatible with the homogeneous essential boundary conditions. Conversely, the space of admissible solutions, $\bfa{V}$, contains all functions from the same Sobolev space that satisfy the prescribed non-homogeneous displacement boundary condition $\bfa{u} = \bfa{u}_0$ on $\Gamma_D$. The definitions are formally stated as:
\begin{subequations}
\begin{align}
\bfa{V}_{\bfzero} &:= \left\{ \bfu \in \left( H^{1}(\Lambda)\right)^2 \colon \; \bfu=\bfzero \quad  \mbox{on} \;\;\Gamma_D\right\}, \label{test_V0}\\
\bfa{V} &:= \left\{ \bfu \in \left( H^{1}(\Lambda)\right)^2 \colon \; \bfu=\bfu_0 \quad \mbox{on} \;\; \Gamma_D\right\}. \label{test_Vu0}
\end{align}
\end{subequations}

The existence of a solution to the nonlinear BVP in \eqref{eq:strong_form_bvp} is established by analyzing its corresponding weak (or variational) formulation, following the framework developed by Beck et al.~\cite{beck2017existence}. The weak formulation is derived by multiplying the equilibrium equation \eqref{eq:equilibrium} by an arbitrary test function $\bfa{w}$ from the space of kinematically admissible variations $\bfa{V_0}$ and integrating over the domain $\Lambda$. Applying the divergence theorem and incorporating the Neumann boundary condition \eqref{eq:neumann_bc} yields the following integral statement:

\begin{cwf} Find the displacement $\bfa{u} \in \bfa{V}$ such that for all test functions $\bfa{w} \in \bfa{V_0}$:
\begin{equation}
\int_{\Lambda} \bfsigma(\bfeps(\bfa{u})) \colon \bfeps(\bfa{w}) \, d\bfa{x} = \int_{\Lambda} \bfa{f} \cdot \bfa{w} \, d\bfa{x} + \int_{\Gamma_N} \bfa{g} \cdot \bfa{w} \, ds. \label{eq:weak_form}
\end{equation}
\end{cwf}
The left-hand side of this equation represents the internal virtual work done by the stresses, while the right-hand side represents the external virtual work done by the body forces and surface tractions. The following theorem, stated without proof, formally asserts the existence of a solution.

\begin{theorem}[Existence of a weak solution]
\label{thm:existence}
Let $\Lambda \subset \mathbb{R}^2$ be a bounded Lipschitz domain, given the body force $\bfa{f} \in (L^2(\Lambda))^2$ and the traction vector $\bfa{g} \in (L^2(\Gamma_N))^2$, and assuming the constitutive law $\bfsigma(\bfeps)$ satisfies the properties of continuity, monotonicity, and coercivity, and that assumptions A1-A3 hold, then there exists at least one solution $\bfa{u} \in \bfa{V}$ to the weak formulation \eqref{eq:weak_form}.
\end{theorem}

This theorem provides the essential theoretical guarantee that the formulated BVP is solvable, paving the way for the development of stable and convergent numerical methods, such as the finite element method, for its approximation. The theoretical underpinning for the existence of a solution to our boundary value problem is provided by the framework established in the work of Beck et al.~\cite{beck2017existence}. The existence theorem presented therein is directly applicable to our current formulation, with only a minor modification in the algebraic form of the Cauchy stress tensor, as specified in equation~\eqref{eq:inverted_hyperelastic}. This congruence holds because our chosen constitutive model preserves the essential mathematical properties, such as strict monotonicity and coercivity, that are central to the theorem's proof. We therefore posit that the material properties ($\mu, \lambda$), the constitutive modeling parameters ($\alpha, \beta$), and the problem data ($\bfa{f}, \bfa{g}, \bfa{u}_0$) all satisfy the necessary conditions stipulated in~\cite{beck2017existence}. With these requirements met, we can confidently assert the existence of a unique solution pair $(\bfu, \bfT)$ residing in the function space $(W^{1,1}(\Lambda))^2 \times \Sym({ \mathbb{L}}^1(\Lambda)^{2 \times 2})$ for the weak formulation previously delineated. The formal proof for our specific strain-limiting model is a direct adaptation of the arguments presented in~\cite{beck2017existence}, providing a rigorous mathematical foundation for the numerical investigations that follow.

\section{Finite element discretization}
\label{sec:fem}

This section is dedicated to the numerical approximation of the nonlinear boundary value problem established in the previous section. Due to the inherent nonlinearity of the constitutive model, obtaining a closed-form analytical solution is generally intractable. Consequently, we employ a conforming Galerkin finite element method to discretize the problem, thereby transforming the continuous variational problem into a tractable system of nonlinear algebraic equations. We begin by formally stating the continuous variational formulation and reviewing the properties of its solution, followed by the construction of the discrete problem that serves as the basis for our numerical algorithm.

\subsection{Variational formulation and solution properties}

The cornerstone of the finite element method is the conversion of the strong form of the BVP, a partial differential equation, into an equivalent integral statement, known as the weak or variational formulation. This is accomplished by applying the method of weighted residuals. The governing equilibrium equation is multiplied by an arbitrary test function $\bfa{v}$ from the space of kinematically admissible variations $\bfa{V}_{\bfzero}$ (defined in \eqref{test_V0}) and integrated over the domain $\Lambda$. Through the application of the divergence theorem (Green's formula) and the incorporation of the Neumann boundary condition \eqref{eq:neumann_bc}, we arrive at the following continuous variational problem.\\

\noindent\textbf{Problem (Continuous Weak Formulation):}
Given the prescribed body forces, surface tractions, and material parameters, find the displacement field $\bfa{u} \in \bfa{V}$ such that the following integral equality holds for all test functions $\bfa{v} \in \bfa{V}_{\bfzero}$:
\begin{equation}
\label{eq:weak_formulation}
a(\bfa{u}; \bfa{v}) = L(\bfa{v}), \quad \forall\, \bfa{v} \in \bfa{V}_{\bfzero},
\end{equation}
where the semilinear form $a(\cdot; \cdot)$ and the linear functional $L(\cdot)$ represent the internal and external virtual work, respectively, and are defined as:
\begin{subequations}
\label{def:A-L}
\begin{align}
a(\bfa{u}; \bfa{v}) &:= \int_{\Lambda} \Psi\left( { \|\mathbb{E}^{1/2} [\bfeps(\bfa{u})]\|}\right) \, \mathbb{E}[\bfeps(\bfa{u})] \colon \bfeps(\bfa{v}) \, d\bfa{x}, \\
L(\bfa{v}) &:= \int_{\Lambda} \bfa{f} \cdot \bfa{v} \, d\bfa{x} + \int_{\Gamma_N} \bfa{g} \cdot \bfa{v} \, ds.
\end{align}
\end{subequations}
It is important to note that the semilinear form $a(\cdot; \cdot)$ is nonlinear with respect to its first argument, $\bfa{u}$, due to the dependence of the function $\Psi$ on the solution's strain field.

Assuming sufficient regularity of the problem data---specifically, for $\bfa{f} \in (L^{2}(\Lambda))^2$, $\bfa{u}_0 \in (H^{1/2}(\Gamma_D))^2$, and $\bfa{g} \in (L^{2}(\Gamma_N))^2$---it can be shown that the variational problem \eqref{eq:weak_formulation} admits a unique weak solution $\bfa{u}$. This solution is known to possess enhanced regularity and resides in the solution space $\mathcal{U}_s := \{ \bfa{v} \in (H^2(\Lambda) \cap W^{1, \infty}(\Lambda))^2 \colon \left. \bfa{v} \right|_{\Gamma_D} = \bfa{u}_0 \}$. Furthermore, the solution satisfies the following \textit{a priori} stability estimate, which confirms the well-posedness of the continuous problem:
\begin{equation}
\| \bfa{u} \|_{H^2(\Lambda)} \leq \widehat{c} \left( \| \bfa{f} \|_{L^2(\Lambda)} + \| \bfa{u}_0 \|_{H^{1/2}(\Gamma_D)} + \| \bfa{g} \|_{L^2(\Gamma_N)} \right), \label{eq:stability_estimate}
\end{equation}
where $\widehat{c} > 0$ is a regularity constant independent of the solution and the problem data.

\begin{remark}
The intrinsic nonlinearity of the constitutive model, represented by the function $\Psi(\cdot)$, necessitates an iterative approach for its numerical solution. In this work, we employ Picard's iterative method (a fixed-point iteration scheme), which is a natural and robust choice for this class of semilinear problems. The convergence of Picard's method is only guaranteed if the initial guess is sufficiently close to the true solution. To ensure a reliable and convergent solution strategy, we first solve the corresponding linear elastic problem (by setting the nonlinearity parameter $\beta = 0$, which simplifies the response function to $\Psi \equiv 1$). The solution of this linear problem then serves as a robust and well-informed initial estimate for the subsequent nonlinear iterations.
\end{remark}

\subsection{Discrete formulation}

To obtain a numerical approximation of the continuous problem, we discretize the computational domain $\overline{\Lambda}$ using a conforming, shape-regular family of partitions, $\{\mathcal{T}_h\}_{h > 0}$, consisting of non-overlapping quadrilateral finite elements $\mathcal{K}$. The mesh parameter $h$ is defined as the maximum diameter of any element in the partition, $h := \max_{\mathcal{K} \in \mathcal{T}_h} \operatorname{diam}(\mathcal{K})$. The set of all boundary edges of the mesh is denoted by $\mathscr{E}_{bd,h}$, which is partitioned into disjoint sets corresponding to the Dirichlet and Neumann boundaries, $\mathscr{E}_{D,h}$ and $\mathscr{E}_{N,h}$, respectively.

We then define a finite-dimensional subspace $\bfa{V}_h \subset (H^1(\Lambda))^2$ for the trial functions and a subspace $\bfa{V}_{h, \, 0} \subset \bfa{V}_{\bfzero}$ for the test functions. The "conforming" nature of this approach is due to the fact that the discrete spaces are subspaces of their continuous counterparts. For this study, we use the space of continuous, vector-valued functions that are piecewise bilinear on each element:
\begin{equation}
\bfa{V}_{h} := \left\{ \bfa{u}_h \in (C^{0}(\overline{\Lambda}))^2 \colon \left. \bfa{u}_h \right|_{\mathcal{K}} \in (\mathbb{Q}_{1}(\mathcal{K}))^2, \forall \mathcal{K} \in \mathcal{T}_h \right\}.
\end{equation}
Here, $\mathbb{Q}_{1}(\mathcal{K})$ is the space of bilinear polynomials on element $\mathcal{K}$, which ensures inter-element continuity and a first-order approximation. The continuous Galerkin method seeks an approximate solution within this discrete space by restricting the continuous weak formulation \eqref{eq:weak_formulation} to the discrete subspaces.

\begin{dwf}
Find the discrete displacement field $\bfa{u}_h \in \bfa{V}_{h}$ (which satisfies the discrete non-homogeneous boundary conditions) such that the following holds for all test functions $\bfa{v}_h \in \bfa{V}_{h,0}$:
\begin{equation}\label{eq:discrete_wf}
a(\bfa{u}_h; \bfa{v}_h) = L(\bfa{v}_h).
\end{equation}
The discrete forms are computed by summing the contributions from each element and boundary edge:
\begin{align}\label{}
a(\bfa{u}_h; \bfa{v}_h) &= \sum_{\mathcal{K} \in \mathcal{T}_h} \int_{\mathcal{K}} \Psi\left( \|\mathbb{E}^{1/2}[\bfeps(\bfa{u}_h)]\| \right) \, \mathbb{E}[\bfeps(\bfa{u}_h)] \colon \bfeps(\bfa{v}_h) \, d\bfa{x}, \label{eq:discrete_a} \\
L(\bfa{v}_h) &= \sum_{\mathcal{K} \in \mathcal{T}_h} \int_{\mathcal{K}} \bfa{f} \cdot \bfa{v}_h \, d\bfa{x} + \sum_{e \in \mathscr{E}_{N,h}} \int_{e} \bfa{g} \cdot \bfa{v}_h \, ds. \label{eq:discrete_L}
\end{align}
\end{dwf}
The well-posedness of this discrete nonlinear system, which guarantees the existence and uniqueness of a discrete solution, is founded on the mathematical properties of the semilinear form $a(\cdot; \cdot)$. One can first establish the \textit{Lipschitz continuity} and \textit{strong monotonicity}, which are essential for proving the existence and uniqueness of the discrete solution and guaranteeing the convergence of iterative solvers used to find the solution \cite{manohar2024hp}.

\section{Numerical results and discussion}
\label{sec:rd}

This section details the transition from the theoretical framework to its practical application through a series of numerical experiments. Our primary objective is to validate the proposed strain-limiting constitutive model and to systematically quantify its impact on the mechanics of a single crack embedded within a transversely isotropic solid. By contrasting the predictions of our nonlinear model with those of classical linear elasticity, we aim to demonstrate how the model's inherent strain-limiting property effectively regularizes the non-physical singularities at the crack tip. This regularization yields a more physically realistic and stable representation of the near-tip stress and strain fields, which is important for accurately predicting material failure. For the numerical approximation, we employ a conventional continuous Galerkin finite element method using a discretization based on bilinear quadrilateral elements. This standard approach is a robust and well-understood methodology, which is perfectly suited for the present goal of demonstrating the fundamental physical differences between the two constitutive models.  All simulations were implemented using the open-source, object-oriented finite element library \texttt{deal.II}~\cite{2023dealii, dealii2019design}, and were performed on structured computational meshes to ensure a high degree of control over the discretization.

The constitutive nonlinearity, the central feature of our model, is resolved using the iterative procedure detailed in Algorithm~\ref{alg:picard}. At each Picard iteration $n$, we assess the convergence of the solution by computing the norm of the residual functional. This residual quantifies the imbalance in the discrete weak formulation and is a direct measure of how close the current approximate solution is to satisfying the equilibrium equation. The iteration is terminated when the norm of this residual drops below a prescribed tolerance. The complete computational procedure for this study is outlined below.

\begin{equation}\label{eq:residual}
\mathcal{R}(\bfa{u}^n_h; \bfa{\varphi}_h) := a(\bfa{u}^n_h; \bfa{\varphi}_h).
\end{equation}

The numerical procedure for solving the nonlinear finite element problem is detailed here. This iterative process is essential for effectively handling the material nonlinearity introduced by the strain-limiting constitutive model. A crucial aspect of this scheme is the monitoring of convergence, which is accomplished by evaluating the norm of the residual functional at each iteration. This residual serves as a direct measure of how closely the current solution satisfies the discrete weak formulation. The iterative loop is terminated when the residual norm falls below a prescribed tolerance, thereby indicating that an equilibrium solution has been reached. The complete computational procedure, including initialization, the iterative loop, and post-processing, is summarized in Algorithm \ref{alg:picard}.

The Picard iteration scheme is a straightforward and effective method for solving the nonlinear problem at hand. As a fixed-point iteration approach, it is well-suited for the class of semilinear problems where the nonlinearity is contained within the constitutive law. However, a key characteristic of this method is its linear convergence rate, which can lead to a significant number of iterations, particularly when the problem exhibits strong nonlinearity. Despite this, the Picard method remains a robust and reliable choice when provided with a sufficiently accurate initial guess, such as the solution from the corresponding linear elastic problem. Here is an algorithm that uses Picard's iterative method to solve the nonlinear problem.\\

\begin{algorithm}[H]
\SetAlgoLined
\KwInput{Finite element mesh $\mathcal{T}_h$; material and model parameters ($\mu, \lambda, \gamma, \alpha, \beta$); max iterations $M_{max}$; tolerance $TOL$; boundary data $\bfa{g}, \bfa{u}_0$.}
\KwOutput{Converged discrete solution $\bfa{u}_h$.}

\tcp{Step 1: Initialization}
Set iteration counter $n \leftarrow 0$. \\
Compute the initial guess $\bfa{u}^0_h$ by solving the linear elastic problem (i.e., setting $\beta=0$ in the weak form).

\tcp{Step 2: Picard's Iteration Loop}
\While{$n < M_{max}$}{
  $n \leftarrow n + 1$. \\
  \tcp{Assemble the stiffness matrix and right-hand side}
  Compute the residual vector $\bfa{R}(\bfa{u}^{n-1}_h)$ using the previous solution. \\
  Compute the stiffness matrix $\bfa{A}(\bfa{u}^n_h)$ by linearizing the weak form.

  \tcp{Step 3: Solve for displacement increment}
  Solve the linear system $\bfa{A}(\bfa{u}^n_h) = \bfa{0}$. \\
  \tcp{Update the solution}
  $\bfa{u}^{n+1}_h \leftarrow \bfa{u}^n_h $.

  \tcp{Step 4: Convergence Check}
  Compute the norm of the displacement increment: $\| \Delta \bfa{u}_h \|$ and the residual. \\
  \If{$\| \Delta \bfa{u}_h \| \leq TOL \quad \& \quad  \bfa{R}(\bfa{u}^{n}_h) \leq TOL $}{
    Break loop.
  }
}

\tcp{Step 5: Output}
Return the converged solution $\bfa{u}_h \leftarrow \bfa{u}^{n+1}_h$. \\
Perform postprocessing (e.g., calculate stress, strain, strain energy).

\caption{Picard's iterative scheme for the nonlinear fracture problem.}
\label{alg:picard}
\end{algorithm}

\subsection{Setup of domain and boundary conditions}

The central numerical experiment of this study is a benchmark problem designed to validate our nonlinear framework: a rectangular plate with a single edge crack subjected to compressive load. Our primary goal is to perform a direct and quantitative comparison of the mechanical fields predicted by our strain-limiting model with the standard predictions of classical linear elastic fracture mechanics. Through this comparison, we aim to demonstrate that the inherent strain-limiting property of our model effectively mitigates the non-physical stress and strain singularities predicted by linear theory, providing a more realistic and bounded depiction of the near-tip fields.

\begin{figure}[H]
  \centering
  \includegraphics[width=0.5\linewidth]{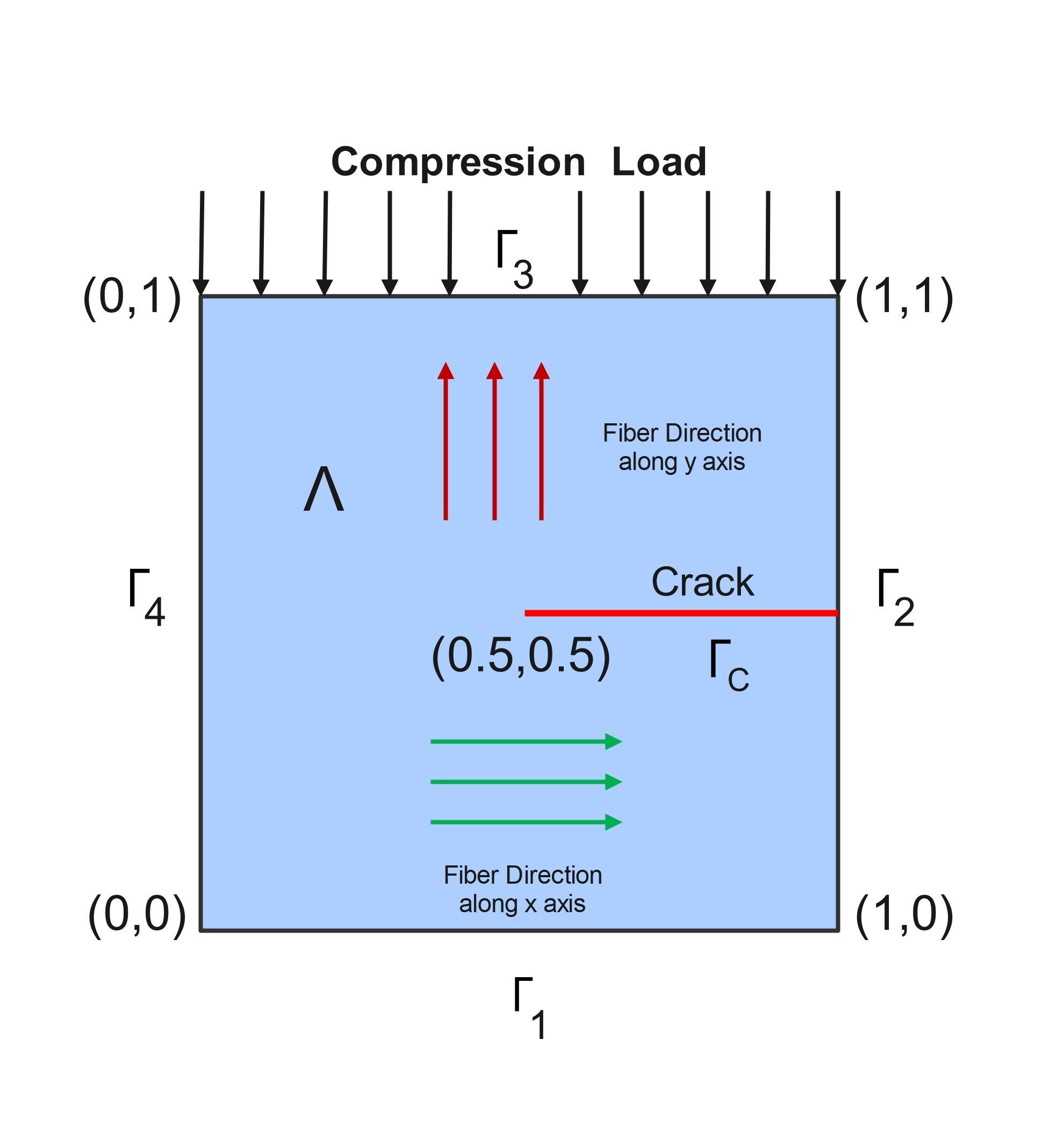}
  \caption{Schematic of the computational domain and boundary indicators.}
  \label{fig:cd}
\end{figure}

The constitutive nonlinearity is handled by the iterative procedure detailed in Algorithm~\ref{alg:picard}. Convergence is assessed at each Picard iteration $n$ by computing the norm of the residual functional, which measures the imbalance in the weak formulation. 
The iterative process is terminated when this residual norm falls below a prescribed tolerance of $TOL = 10^{-4}$, with a maximum of $M_{max} = 200$ iterations, which was sufficient for convergence in all cases presented. The geometric setup and boundary conditions for this benchmark problem are illustrated in Figure~\ref{fig:cd}, where the computational domain is a rectangular plate with a horizontal crack extending from the left edge along (parallel to the $x$-axis) along a line $y=0.5, \;\; 0.5 \leq x \leq 1$. 

The prescribed loads and constraints are as follows:
\begin{itemize}
 \item \textbf{Top boundary ($\Gamma_{3}$):} A compressive load ($u_2 = - c$) is applied.
 \item \textbf{Right boundary ($\Gamma_{2}$):} This edge is traction-free.
\item \textbf{Left boundary ($\Gamma_{4}$):} A traction-free boundary condition is considered on this edge.
\item \textbf{Bottom boundary ($\Gamma_{1}$):} Vertical movement is restricted ($u_2 = 0$), enforcing a roller kind of boundary condition ($u_1$ is not constrained). 
    \item \textbf{Crack faces ($\Gamma_{c}$):} The crack faces are free of traction. 
\end{itemize}

A defining characteristic of transversely isotropic materials is their directionally dependent mechanical response, which is fundamentally governed by their internal microstructure. In our mathematical framework, this anisotropy is explicitly captured through the orientation of the structural tensor, $\bfa{M}$, within the fourth-order elasticity tensor, $\mathbb{E}$. The alignment of this structural tensor dictates the material's principal axis of stiffness. To systematically investigate the influence of this anisotropy on crack-tip phenomena, our study explores two distinct scenarios for the material's fiber orientation. This approach is physically motivated, as it is well-established that the axis of material symmetry in materials such as fiber-reinforced composites, wood, and bone is determined by the alignment of their constituent fibers or other microstructural features. The orientation of these internal structures plays a crucial role in determining the material's overall mechanical behavior. It fundamentally governs how stresses are distributed and concentrated under applied loads, particularly in the critical region surrounding a crack tip. By analyzing these different configurations, we aim to gain critical insight into how material anisotropy can influence stress shielding, strain localization, and ultimately, the fracture resistance of the solid.

\subsection{Isotropic solid under compressive load}
When an isotropic solid is subjected to a compressive load, its dimensions change predictably. Because the material is isotropic, its properties are identical no matter which direction you test them in. This means that the relationship between stress and strain is the same regardless of how the load is applied. When compressed, the solid uniformly shrinks along the axis of the applied force and expands in all directions perpendicular to it. This predictable deformation is a defining characteristic of isotropic materials under load. 

\begin{figure}[H]
    \centering 
    \begin{subfigure}[b]{0.48\textwidth}
        \includegraphics[width=\linewidth]{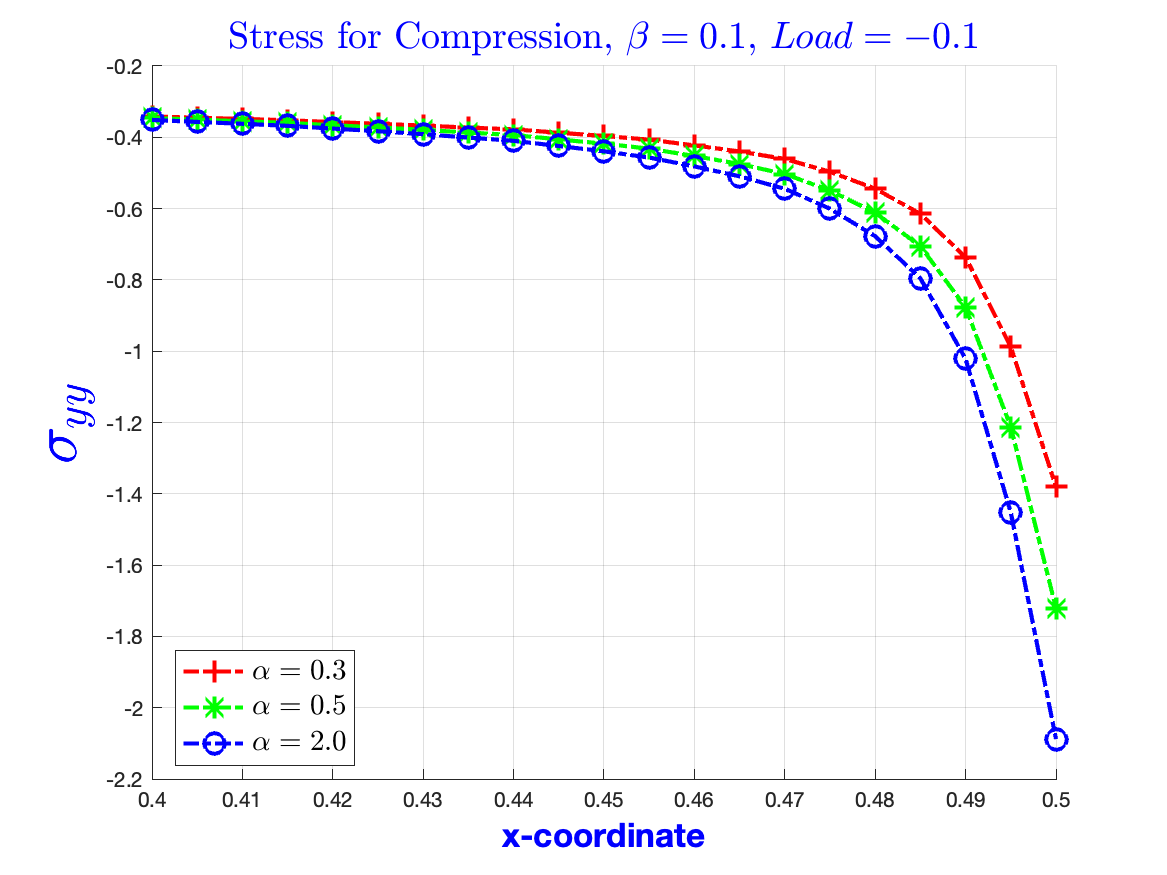}
        \caption{Effect of $\alpha$ on $\bfsigma_{yy}$ }
        \label{fig:fig1}
    \end{subfigure}
    \hfill 
    \begin{subfigure}[b]{0.48\textwidth}
        \includegraphics[width=\linewidth]{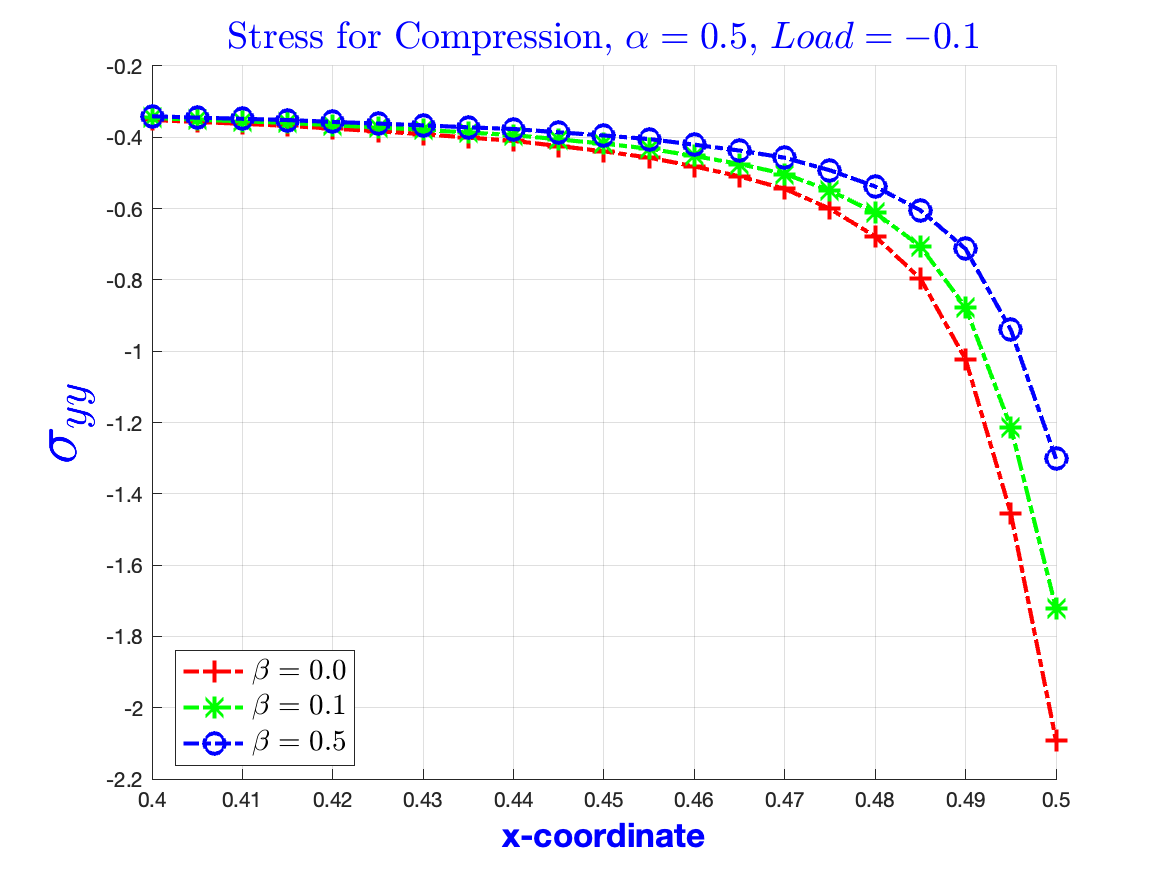}
        \caption{Effect of $\beta$ on $\bfsigma_{yy}$}
        \label{fig:fig2}
    \end{subfigure}
    \vspace{0.5cm} 
    \begin{subfigure}[b]{0.48\textwidth}
        \includegraphics[width=\linewidth]{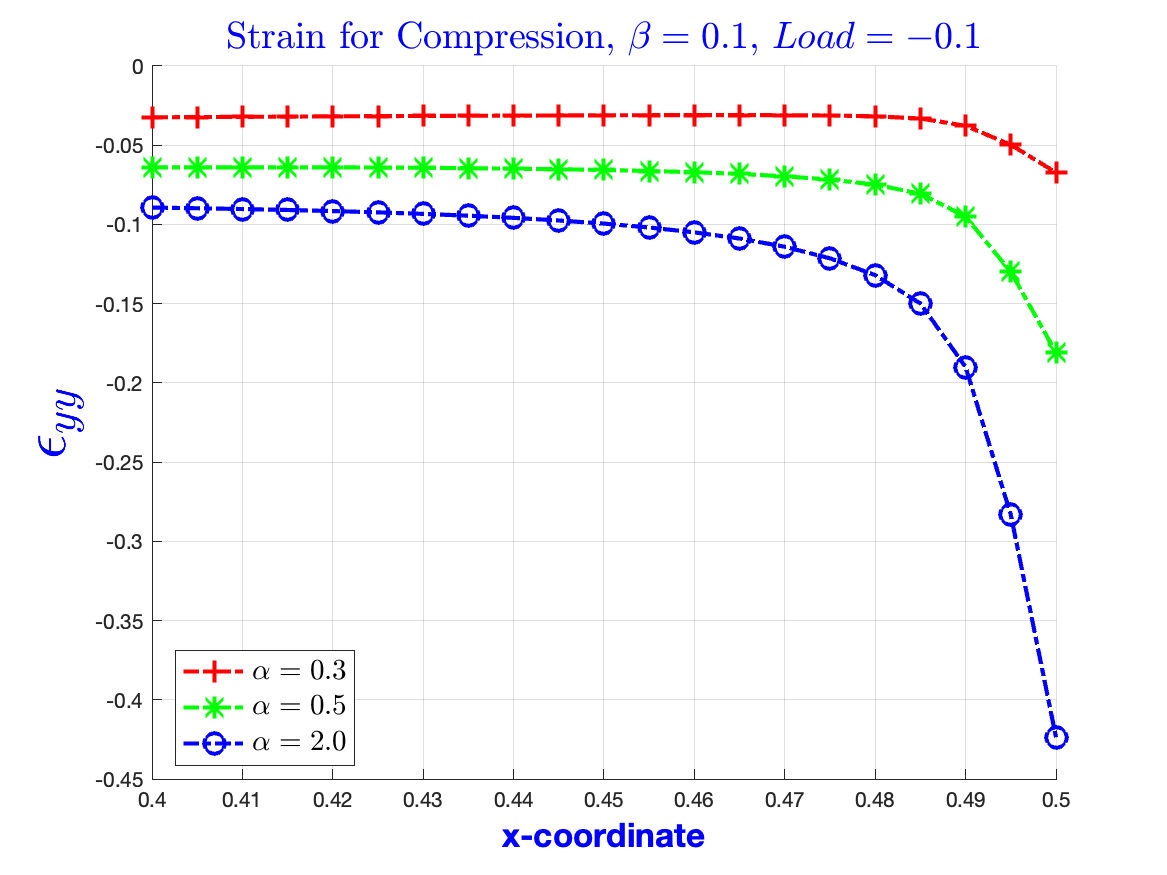}
        \caption{Effect of $\alpha$ on $\bfa{\epsilon}_{yy}$}
        \label{fig:fig3}
    \end{subfigure}
    \hfill 
    \begin{subfigure}[b]{0.48\textwidth}
        \includegraphics[width=\linewidth]{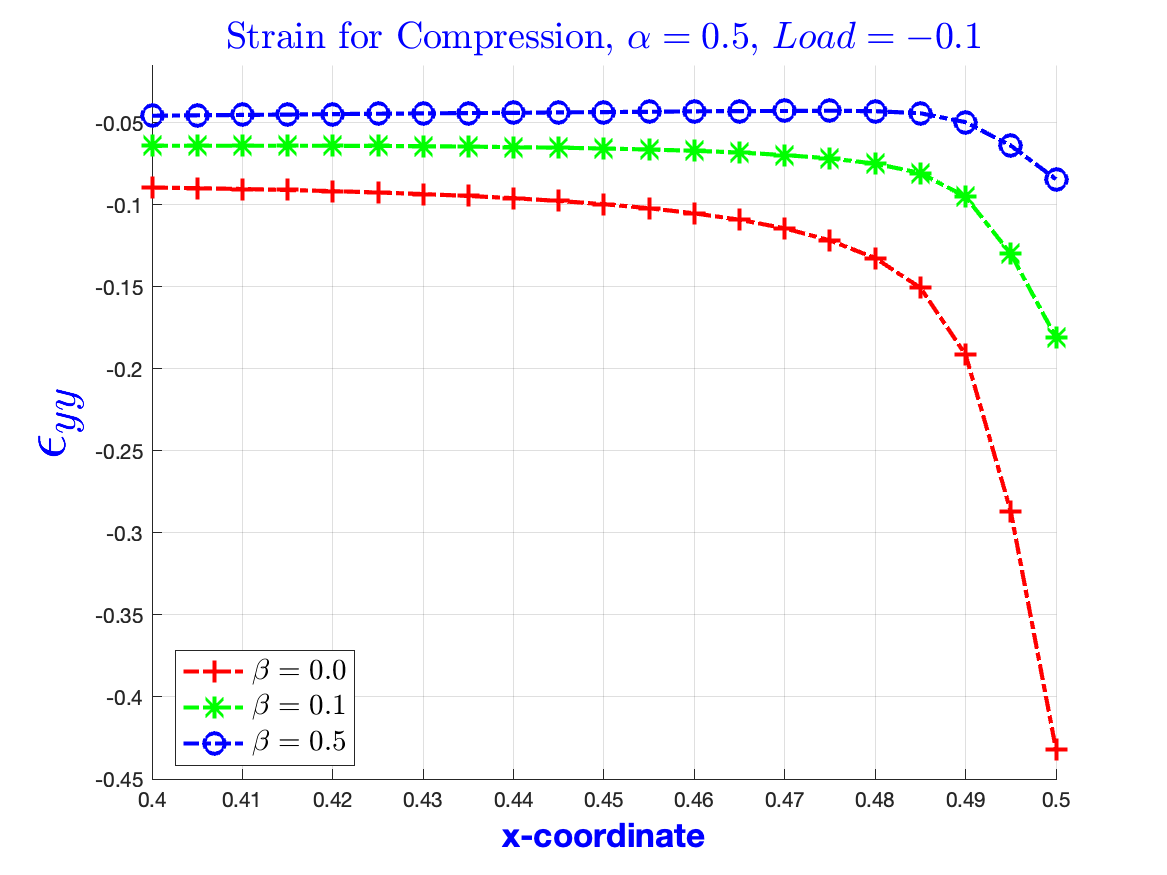}
        \caption{Effect of $\beta$ on $\bfa{\epsilon}_{yy}$}
        \label{fig:fig4}
    \end{subfigure}
    \vspace{0.5cm} 
    \begin{subfigure}[b]{0.48\textwidth}
        \includegraphics[width=\linewidth]{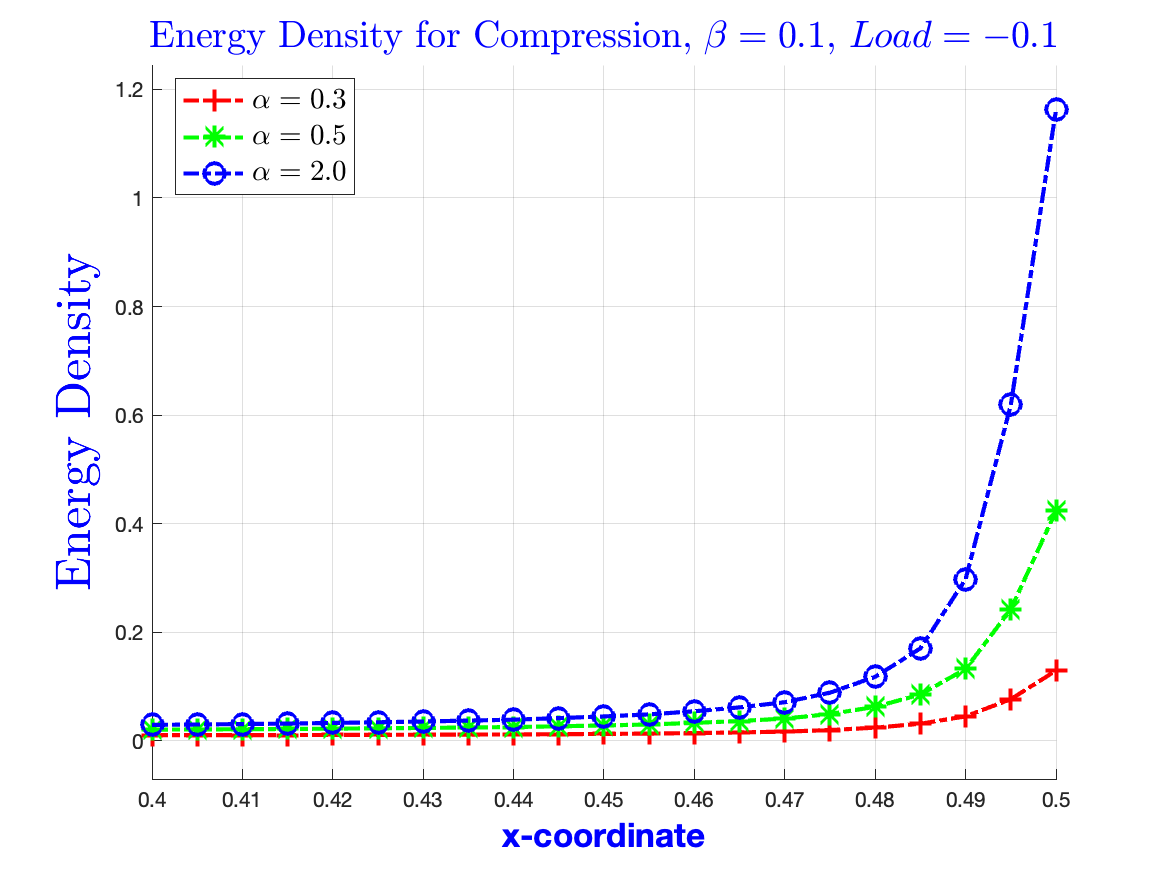}
        \caption{Effect of $\alpha$ on strain-energy density}
        \label{fig:fig5}
    \end{subfigure}
    \hfill 
    \begin{subfigure}[b]{0.48\textwidth}
        \includegraphics[width=\linewidth]{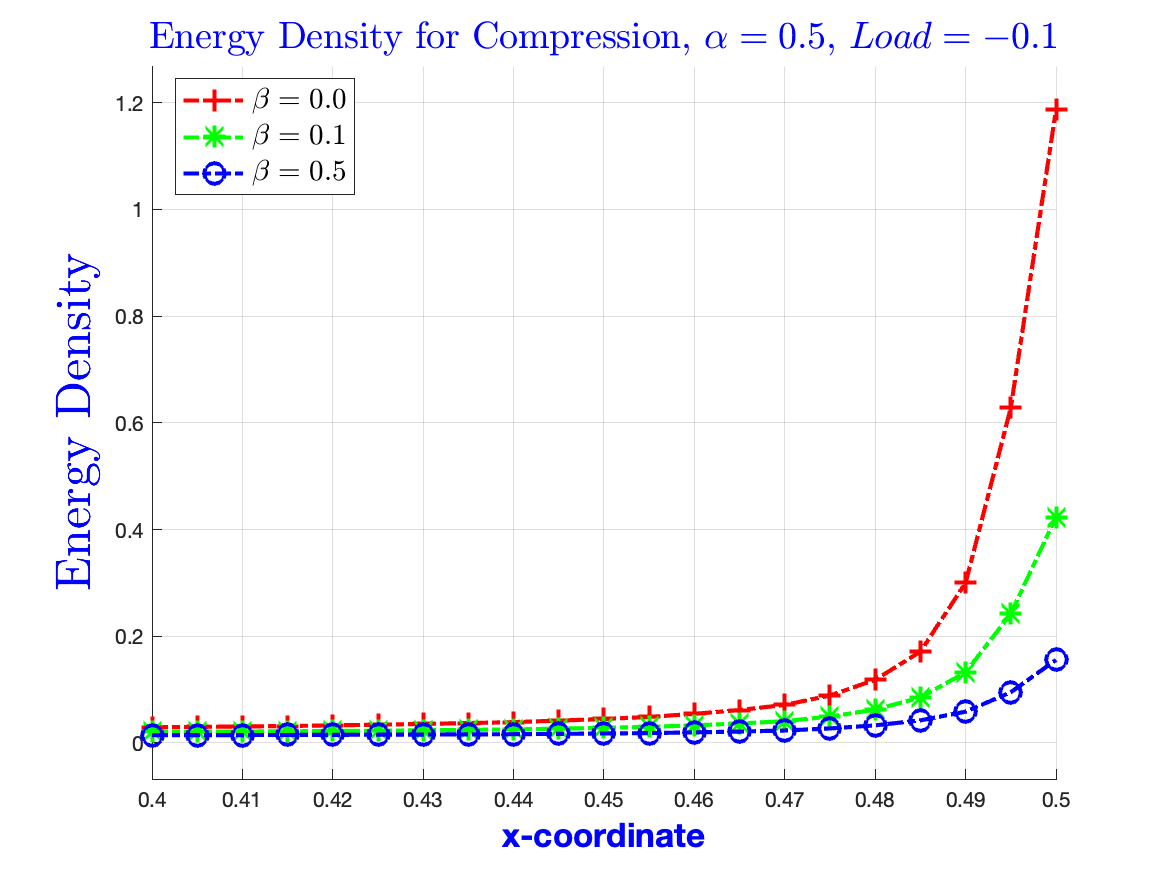}
        \caption{Effect of $\beta$ on strain-energy density}
        \label{fig:fig6}
    \end{subfigure}
\caption{This two-panel figure illustrates the effect of the modeling parameters $\alpha$ and $\beta$ on the crack-tip fields. The left panel shows how $\alpha$ influences key crack-tip quantities like stress, strain, and strain-energy density, while the right panel presents the effect of $\beta$ on the same fields.}
    \label{fig:stress_strain_iso}
\end{figure}

Figure~\ref{fig:stress_strain_iso} illustrates the distinct influence of the two primary modeling parameters, $\alpha$ and $\beta$, on the crack-tip fields under a fixed top compressive load. The computations were performed for various values of these parameters to isolate their individual effects on the crack-tip stress, strain, and strain-energy density. The left panel of the figure specifically depicts the effect of $\alpha$. A key finding is that varying $\alpha$ produces a similar concentration of compressive stress, strain, and strain-energy density around the crack tip. This suggests that the parameter $\alpha$ is instrumental in recovering the results of classical linear elasticity from the model. In a practical engineering sense, a particular value of $\alpha$ can be calibrated to match established material behavior, serving as a baseline for the model's predictive capabilities. In contrast, the right panel highlights the significant and distinct impact of the parameter $\beta$. As the value of $\beta$ increases, a slight decrease in near-tip compressive stresses is observed. However, this is paradoxically accompanied by a reduction in both crack-tip strain and the corresponding strain-energy density. This behavior is of critical engineering importance, as it indicates a shift from a ductile, high-strain state to a more brittle, high-stress state. For a material under compressive loading, an increase in local stress without a corresponding increase in strain and energy absorption suggests a diminished capacity to deform plastically before failure. This makes the material more susceptible to sudden, brittle failure modes, such as crushing or axial splitting, where localized high stresses cannot be relieved by plastic deformation. The parameter $\beta$ thus provides a mechanism to model material behavior where the crack-tip fields are highly localized and do not follow the simple stress-strain relationship of classical elasticity. For instance, this could be used to simulate the behavior of brittle materials or materials under high-rate loading where the strain response is decoupled from the stress concentration. The ability to independently control these two parameters allows for a more nuanced and accurate representation of complex material failure mechanisms.
 
\begin{figure}[H]
    \centering
    \begin{subfigure}{0.3\linewidth}
        \centering
        \includegraphics[width=\linewidth]{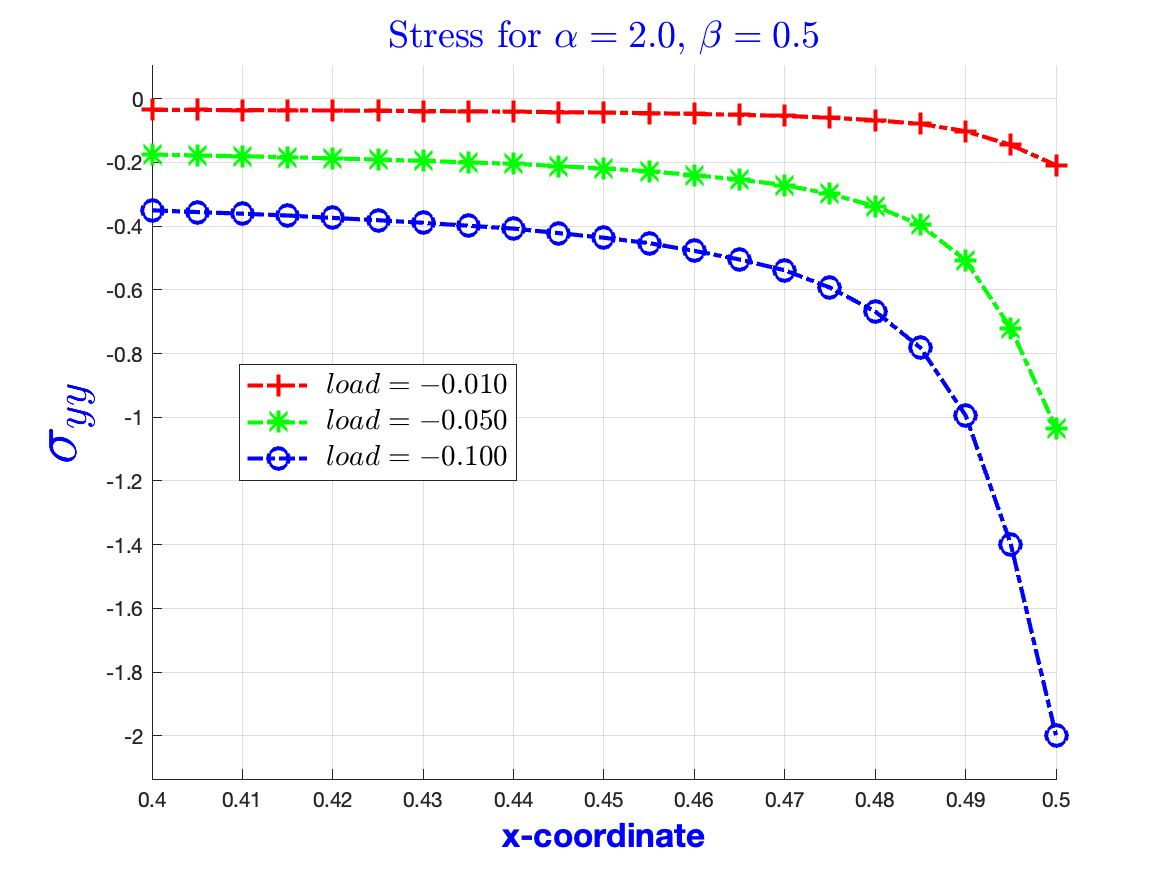}
        \caption{Stress }
        \label{fig:strain_beta_iso}
    \end{subfigure}
    \hfill
    \begin{subfigure}{0.3\linewidth}
        \centering
        \includegraphics[width=\linewidth]{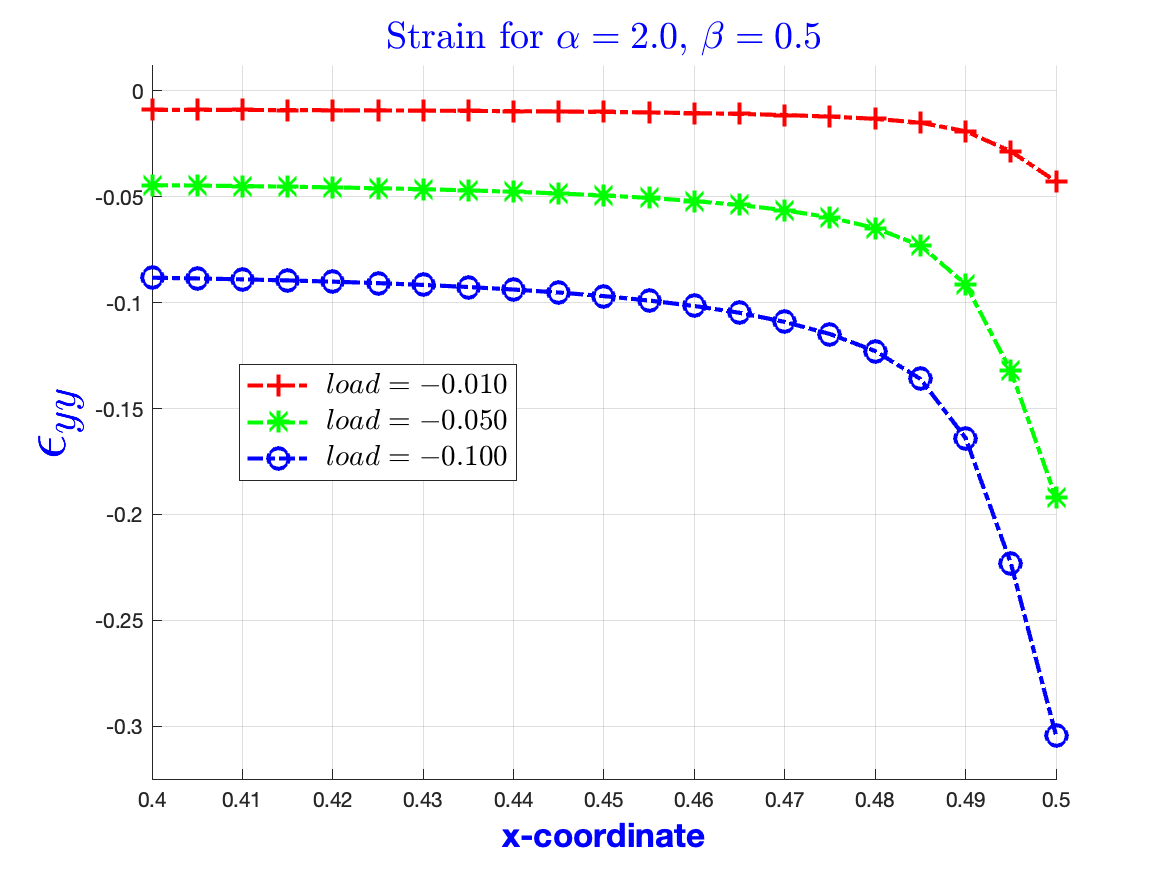}
        \caption{Strain }
        \label{fig:strain_alpha_iso}
    \end{subfigure}
    \hfill
    \begin{subfigure}{0.3\linewidth}
        \centering
        \includegraphics[width=\linewidth]{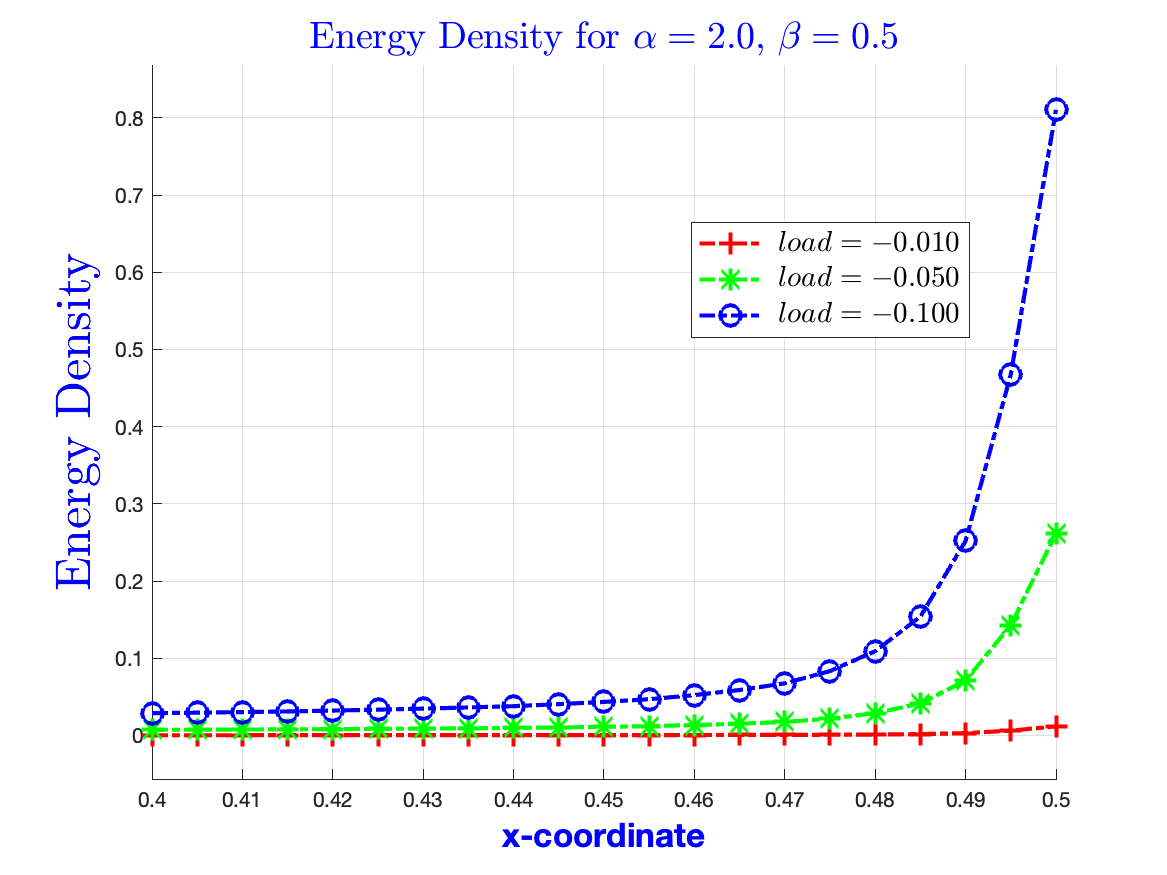}
        \caption{Strain energy density  }
        \label{fig:strain_sigma_iso}
    \end{subfigure}
    \caption{Crack-tip fields for various values of applied top load  }
    \label{fig:SE_load_iso}
\end{figure}

Figure~\ref{fig:SE_load_iso} illustrates the direct effect of the applied top compressive load on the crack-tip fields, with the modeling parameters fixed at $\alpha = 2.0$ and $\beta = 0.5$. The results clearly show a proportional relationship between the magnitude of the applied load and the resulting crack-tip stress, strain, and strain-energy density. As the top compressive load increases, all three of these crack-tip quantities rise significantly. From an engineering perspective, this is a critical observation. The intensification of crack-tip fields with increasing load directly translates to a heightened risk of material failure. Specifically, the observed increase in compressive stress at the crack tip means that the material is being pushed closer to its ultimate compressive strength. Even for materials that are robust under compression, the presence of a crack or flaw creates a stress concentration that can locally exceed the material's failure criteria. This increased susceptibility to failure is particularly relevant for brittle materials, where fracture can initiate from such a stress concentration without significant plastic deformation. In practical applications, this data underscores the importance of a conservative design approach, where the maximum operating load must be kept well below the threshold that would induce critical crack-tip stress levels, thereby preventing catastrophic compressive failure modes such as buckling, crushing, or splitting of the component. The figure serves as a good validation of the model's ability to capture this fundamental engineering principle.

\begin{figure}[H]
    \centering
    % First row
    \begin{subfigure}{0.45\linewidth}
        \centering
        \includegraphics[width=\linewidth]{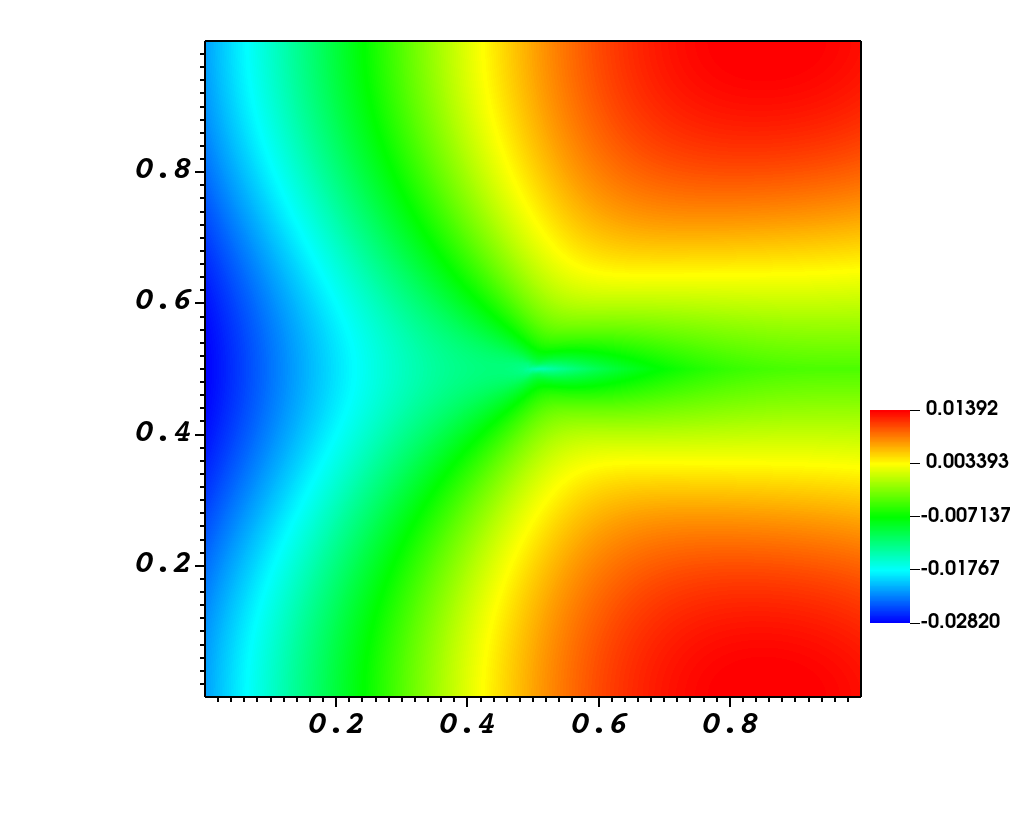}
        \caption{$u_1$}
        \label{fig:x_displacement}
    \end{subfigure}
    \hfill
    \begin{subfigure}{0.45\linewidth}
        \centering
        \includegraphics[width=\linewidth]{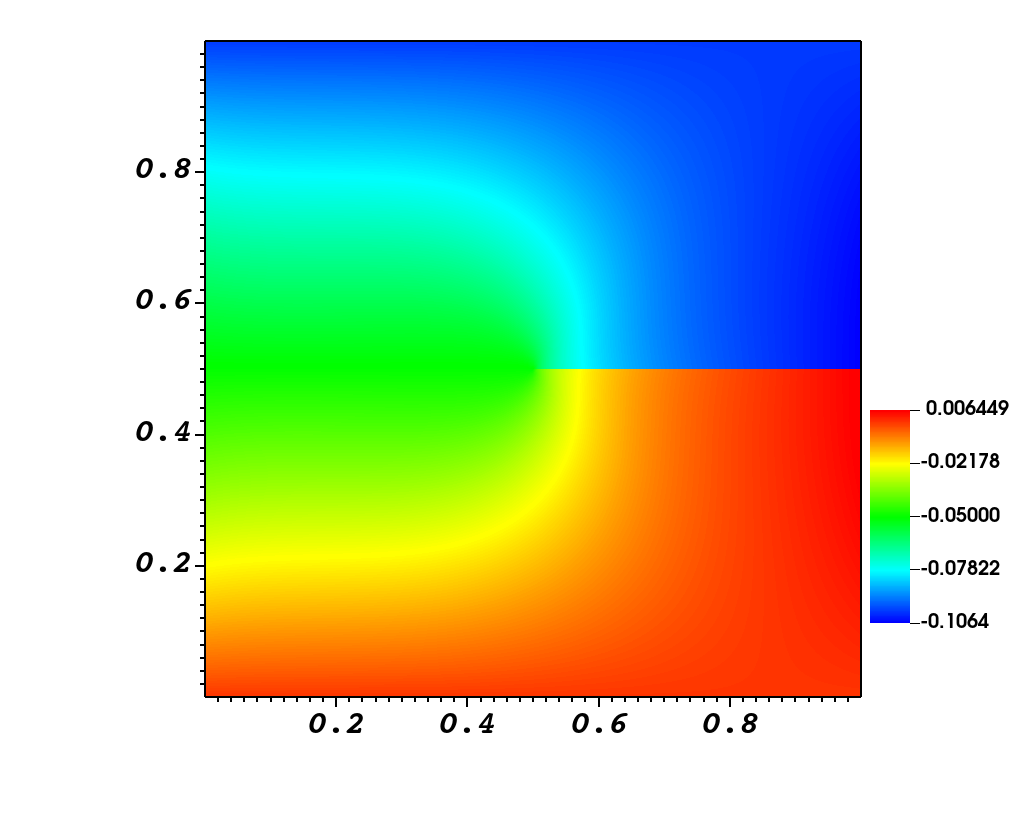}
        \caption{$u_2$ }
        \label{fig:y_displacement}
    \end{subfigure}
     \caption{Displacement plots for  compression load  in isotropic case.}
    \label{fig:displacement_Isotropy}
\end{figure}

Figure~\ref{fig:displacement_Isotropy} presents the displacement components, $u_1$ (horizontal) and $u_2$ (vertical), for a specific set of computational parameters: $\alpha = 2.0$, $\beta = 0.1$, and an applied top compressive load $u_2 = -0.1$. The visual representation of the displacement fields provides crucial insight into the mechanical behavior of the material under load. The plot of the horizontal displacement, $u_1$, clearly illustrates the effect of the applied boundary conditions. The "rolling" boundary condition, where horizontal movement is free but vertical displacement ($u_2$) is fixed at the bottom boundary, is evident. This condition is frequently used in engineering simulations to prevent rigid body motion while allowing the material to deform naturally in the horizontal direction. This is particularly important for accurately capturing the stress distribution near the crack. The plot for the vertical displacement, $u_2$, provides an equally important observation. The displacement field reveals a hint of crack face sliding, which is characteristic of shear-dominated loading conditions. While the primary load is compressive, the geometry of the crack and the boundary conditions induce a localized shear component, causing the crack faces to slide relative to one another. This observation is not only a typical finding for this type of problem but also a critical factor in fracture mechanics. In real-world engineering applications, crack face sliding can lead to complex failure modes, including crack propagation under mixed-mode loading. Understanding this behavior is essential for predicting the fatigue life and ultimate strength of components, particularly in materials where shear failure is a concern. The displacement fields, therefore, serve as a validation of the model's ability to accurately capture both the applied boundary conditions and the resultant complex deformation patterns.

\subsection{Orthotropic solid with fiber orientation along $x$-axis}
In this case, the material is treated as transversely isotropic, meaning its properties are uniform in a plane but different along the axis perpendicular to that plane. Here, the fibers are oriented along the x-axis, which also serves as the axis of symmetry and is parallel to the plane of the crack. This specific fiber orientation is represented by the structural tensor $\bfa{M} = \bfa{e}_1 \otimes \bfa{e}_1$, where $\bfa{e}_1$ is the unit vector along the x-axis. From an engineering standpoint, this setup is crucial for modeling composite materials, such as carbon fiber reinforced polymers, where the strength and stiffness are highly dependent on fiber alignment. By setting the fibers parallel to the crack plane, the model can investigate how a crack behaves when it is aligned with the material's stiffest direction.

\begin{figure}[H]
    \centering 
    \begin{subfigure}[b]{0.48\textwidth}
        \includegraphics[width=\linewidth]{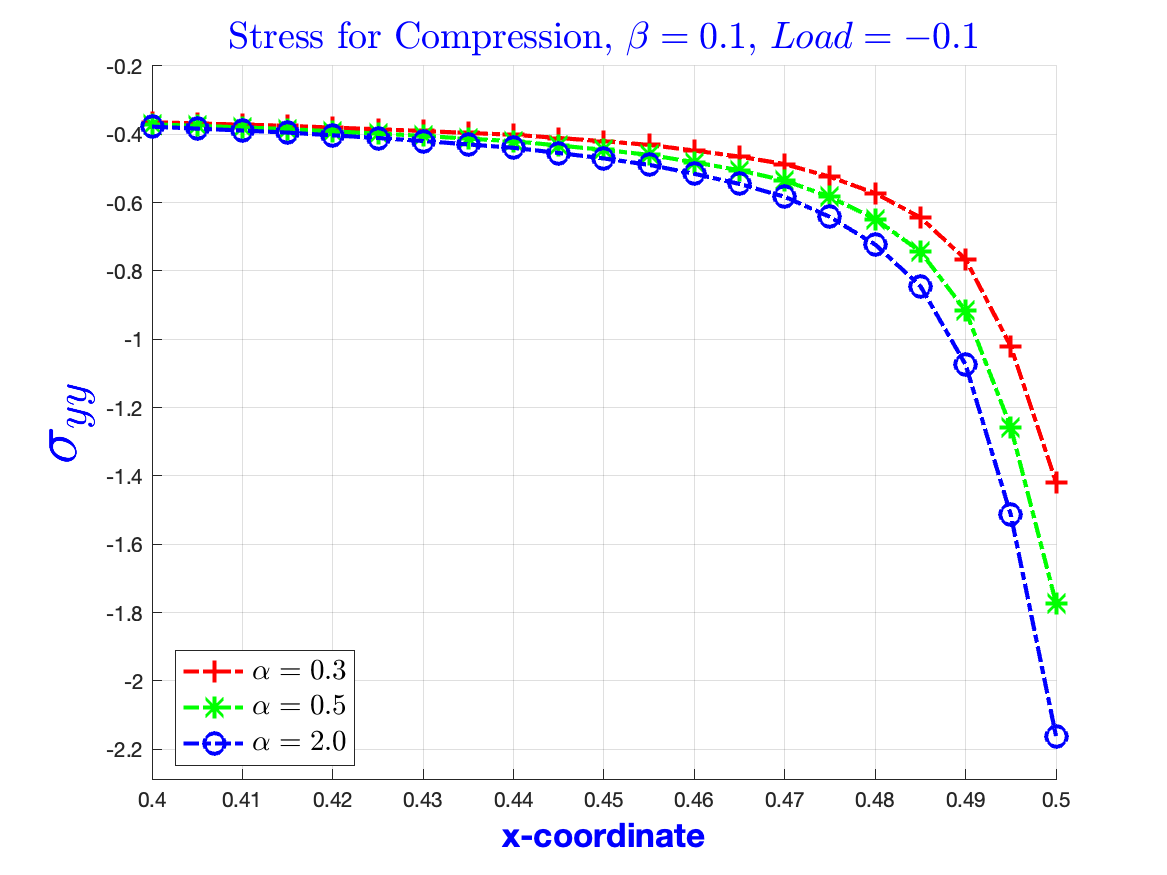}
        \caption{Effect of $\alpha$ on $\bfsigma_{yy}$ }
        \label{fig:fig1}
    \end{subfigure}
    \hfill 
    \begin{subfigure}[b]{0.48\textwidth}
        \includegraphics[width=\linewidth]{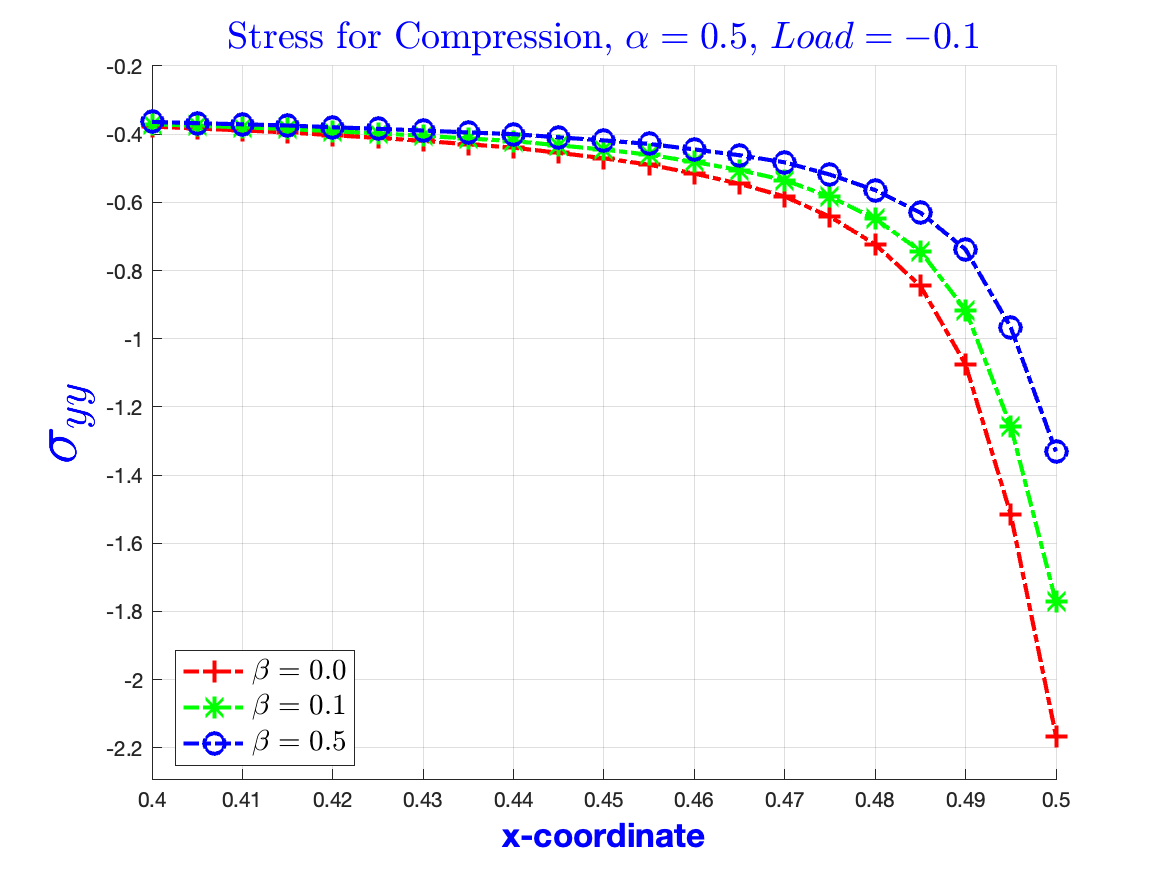}
        \caption{Effect of $\beta$ on $\bfsigma_{yy}$}
        \label{fig:fig2}
    \end{subfigure}
    \vspace{0.5cm} 
    \begin{subfigure}[b]{0.48\textwidth}
        \includegraphics[width=\linewidth]{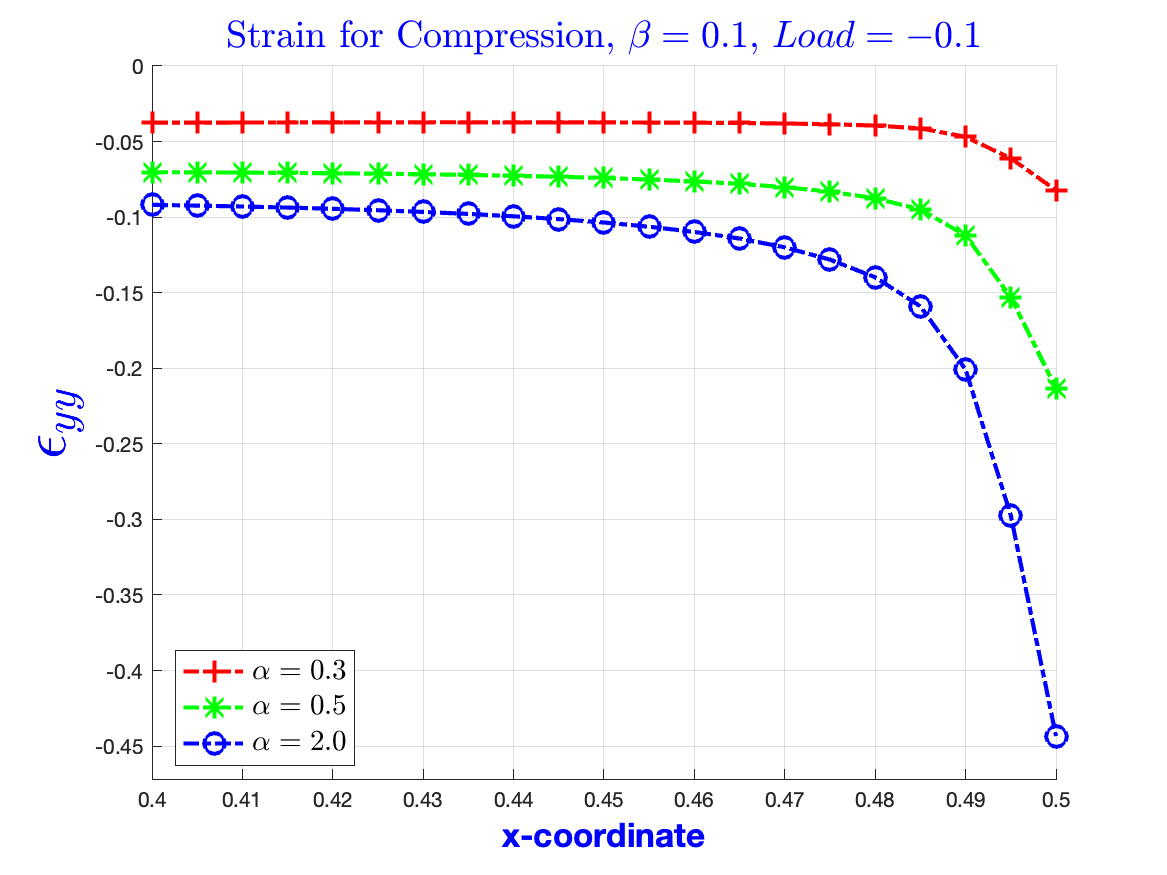}
        \caption{Effect of $\alpha$ on $\bfa{\epsilon}_{yy}$}
        \label{fig:fig3}
    \end{subfigure}
    \hfill 
    \begin{subfigure}[b]{0.48\textwidth}
        \includegraphics[width=\linewidth]{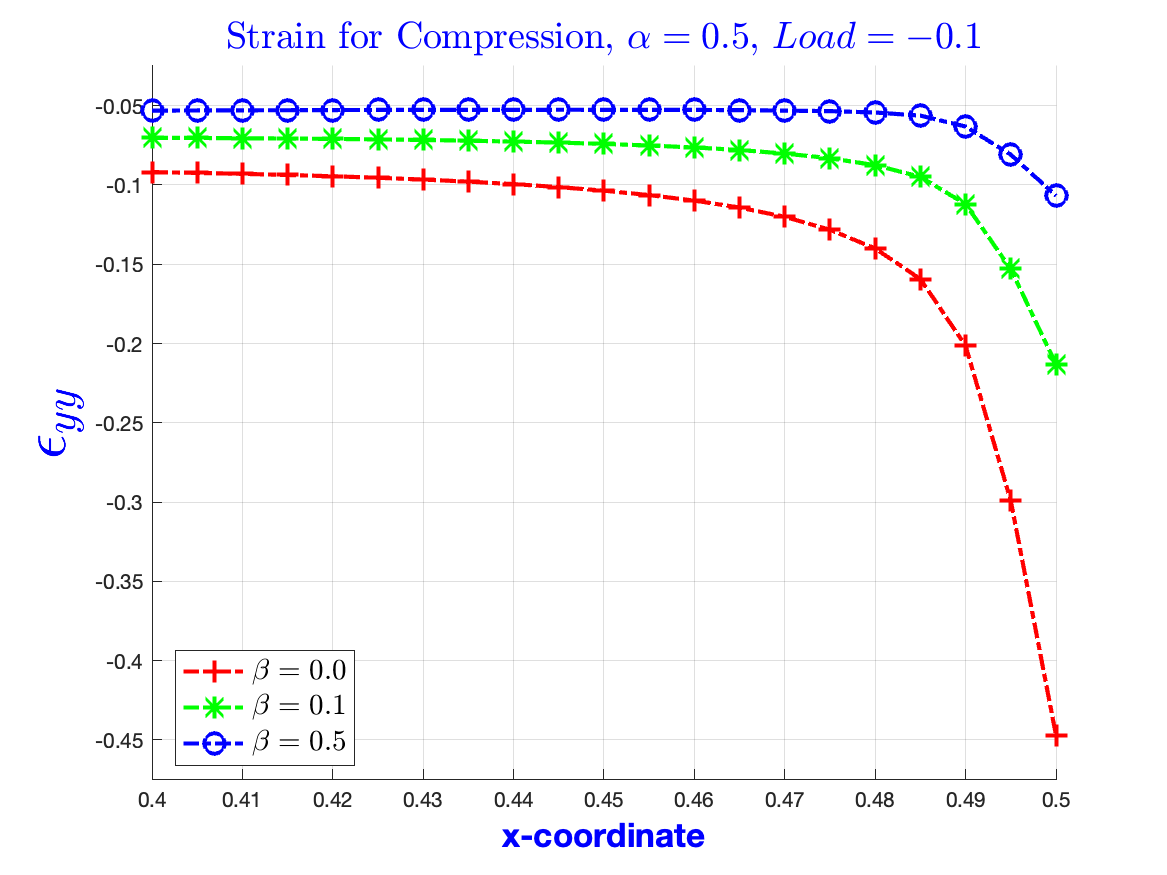}
        \caption{Effect of $\beta$ on $\bfa{\epsilon}_{yy}$}
        \label{fig:fig4}
    \end{subfigure}
    \vspace{0.5cm} 
    \begin{subfigure}[b]{0.48\textwidth}
        \includegraphics[width=\linewidth]{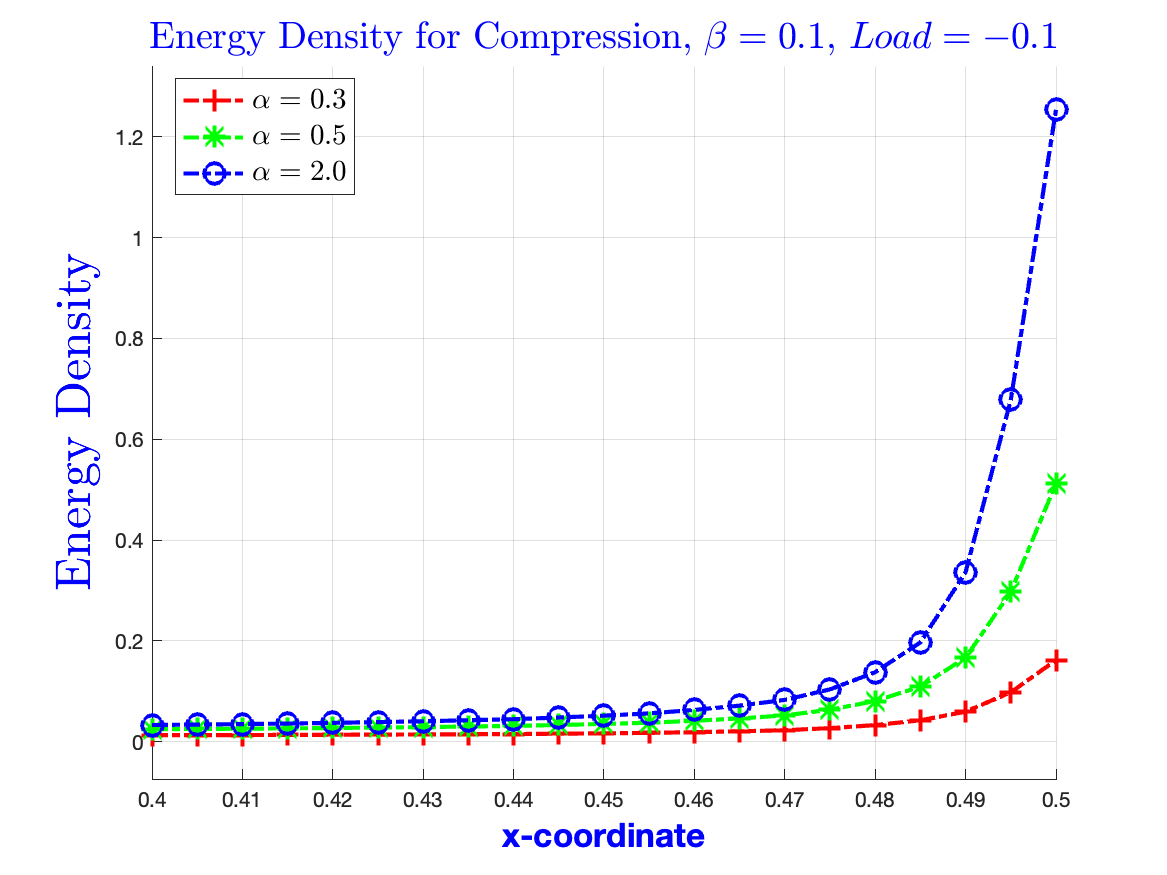}
        \caption{Effect of $\alpha$ on strain-energy density}
        \label{fig:fig5}
    \end{subfigure}
    \hfill 
    \begin{subfigure}[b]{0.48\textwidth}
        \includegraphics[width=\linewidth]{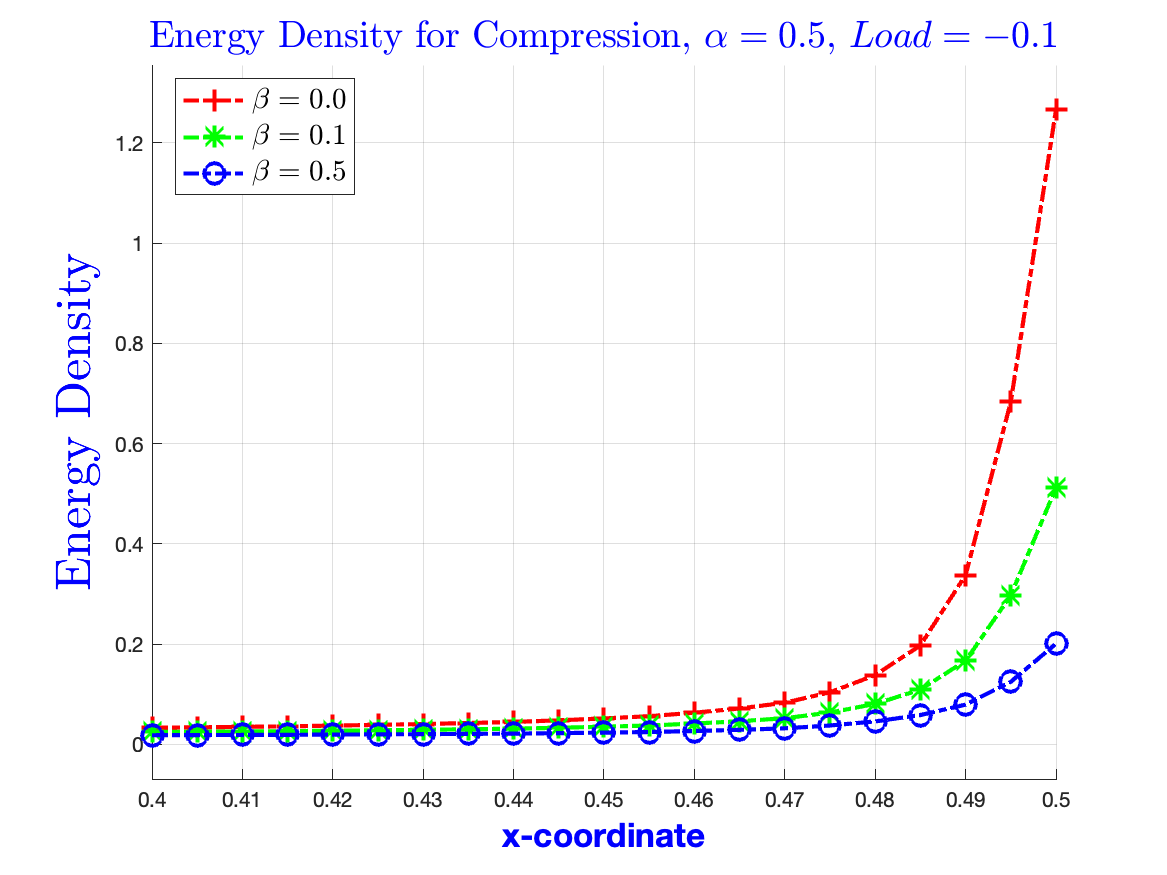}
        \caption{Effect of $\beta$ on strain-energy density}
        \label{fig:fig6}
    \end{subfigure}
\caption{This two-panel figure comprehensively illustrates the influence of two key modeling parameters, $\alpha$ and $\beta$, on the crack-tip fields within an orthotropic solid. The material is specifically configured with its fibers aligned parallel to the crack.}
    \label{fig:stress_strain_M1}
\end{figure}

Figure~\ref{fig:stress_strain_M1} provides a comprehensive look at how the two key modeling parameters, $\alpha$ and $\beta$, influence the crack-tip fields within a strain-limiting material. The material is configured with an edge crack, and its internal fibers are intentionally aligned parallel to the crack. This specific orthotropic setup is crucial for simulating the mechanical behavior of real-world composite materials, such as those found in aerospace or automotive components, where fiber orientation dictates strength and stiffness. The computations were performed with a fixed top compressive load, allowing us to isolate and understand the individual effects of $\alpha$ and $\beta$. The left panel of the figure focuses on the impact of $\alpha$. The results show that varying $\alpha$ produces a proportional and similar concentration of compressive stress, strain, and strain-energy density around the crack tip. From an engineering standpoint, this indicates that $\alpha$ can be used to recover and validate the model's behavior against the well-established principles of classical linear elasticity. By adjusting $\alpha$, engineers can calibrate the model to match the fundamental elastic response of a material, providing a crucial baseline for more complex analyses. Conversely, the right panel highlights the distinct and more complex influence of $\beta$. As the value of $\beta$ increases, the model predicts a paradoxical behavior: a higher concentration of compressive stress at the crack tip is accompanied by a reduction in both compressive strain and strain-energy density. This is a critical finding for engineering applications. It suggests that $\beta$ governs the material's susceptibility to brittle failure under compression. When a material exhibits high stress but low strain and energy absorption, it has a diminished capacity to deform plastically. This behavior indicates a shift towards a more brittle failure mode, such as sudden crushing or splitting, where high localized stresses cannot be dissipated through deformation. The parameter $\beta$ therefore offers a powerful tool for simulating materials where the stress response is highly localized and decoupled from the strain response, enabling more accurate predictions for brittle composites or materials under extreme loading conditions.

\subsection{Orthotropic solid with fiber orientation orthogonal to $x$-axis}
In this case, the material's behavior is defined by transverse isotropy, where its mechanical properties are consistent in one plane but different along the axis perpendicular to that plane. A key aspect of this setup is the deliberate orientation of the material's fibers. Here, the fibers are aligned parallel to the $y$-axis, which is also the axis of symmetry and, crucially, is oriented perpendicular to the crack. This fiber orientation is represented by the structural tensor $\bfa{M} = \bfa{e}_2 \otimes \bfa{e}_2$, where $\bfa{e}_2$ is the unit vector parallel to the $y$-axis. From an engineering perspective, this configuration is vital for understanding how materials like fiber-reinforced composites behave when a flaw, or crack, is oriented perpendicular to their stiffest direction. This setup simulates a worst-case scenario where the crack is most likely to propagate through the weaker matrix material rather than being arrested by the stronger fibers. By modeling this specific orientation, engineers can gain critical insights into the material's fracture toughness and its susceptibility to failure under different loading conditions. The model allows for a detailed analysis of how the stress and strain fields concentrate around the crack tip when the primary load-bearing fibers are unable to bridge the crack, providing essential data for designing safer and more reliable components.

\begin{figure}[H]
    \centering 
    \begin{subfigure}[b]{0.48\textwidth}
        \includegraphics[width=\linewidth]{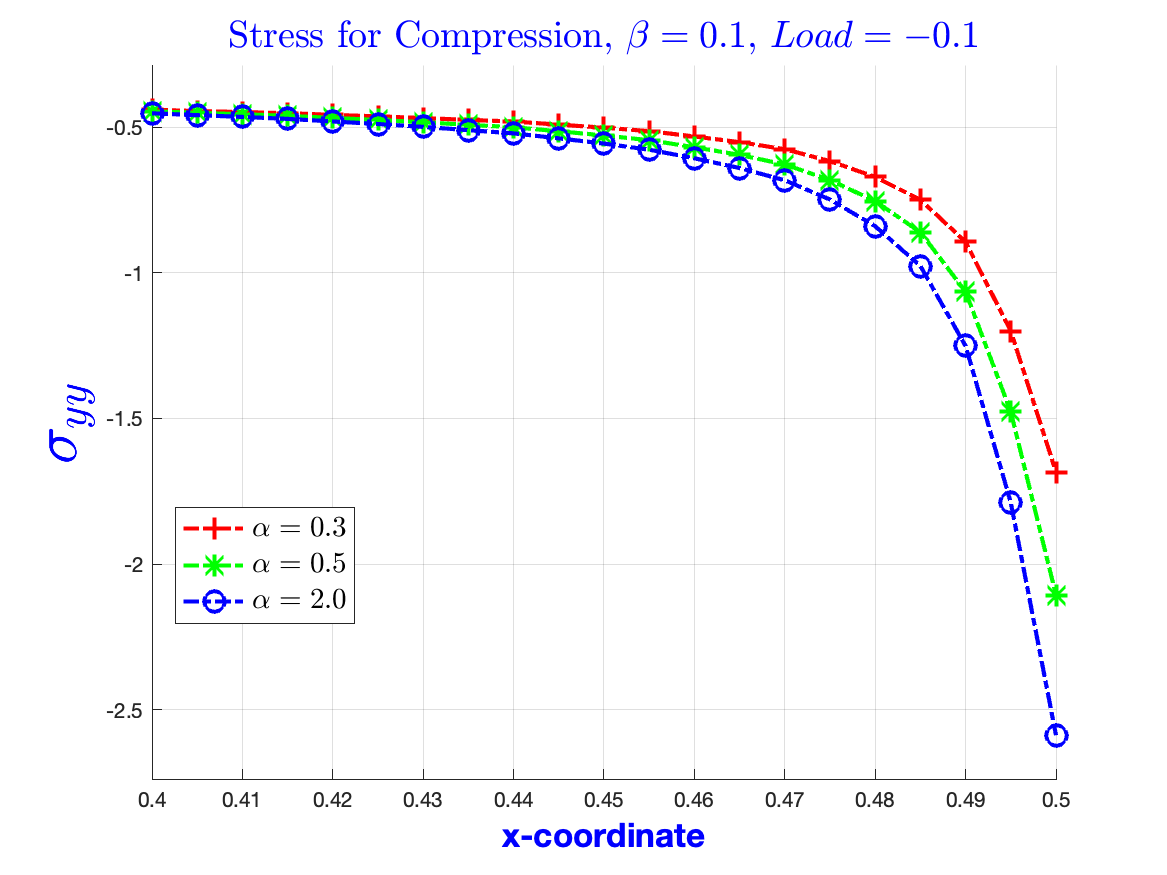}
        \caption{Effect of $\alpha$ on $\bfsigma_{yy}$ }
        \label{fig:fig1}
    \end{subfigure}
    \hfill 
    \begin{subfigure}[b]{0.48\textwidth}
        \includegraphics[width=\linewidth]{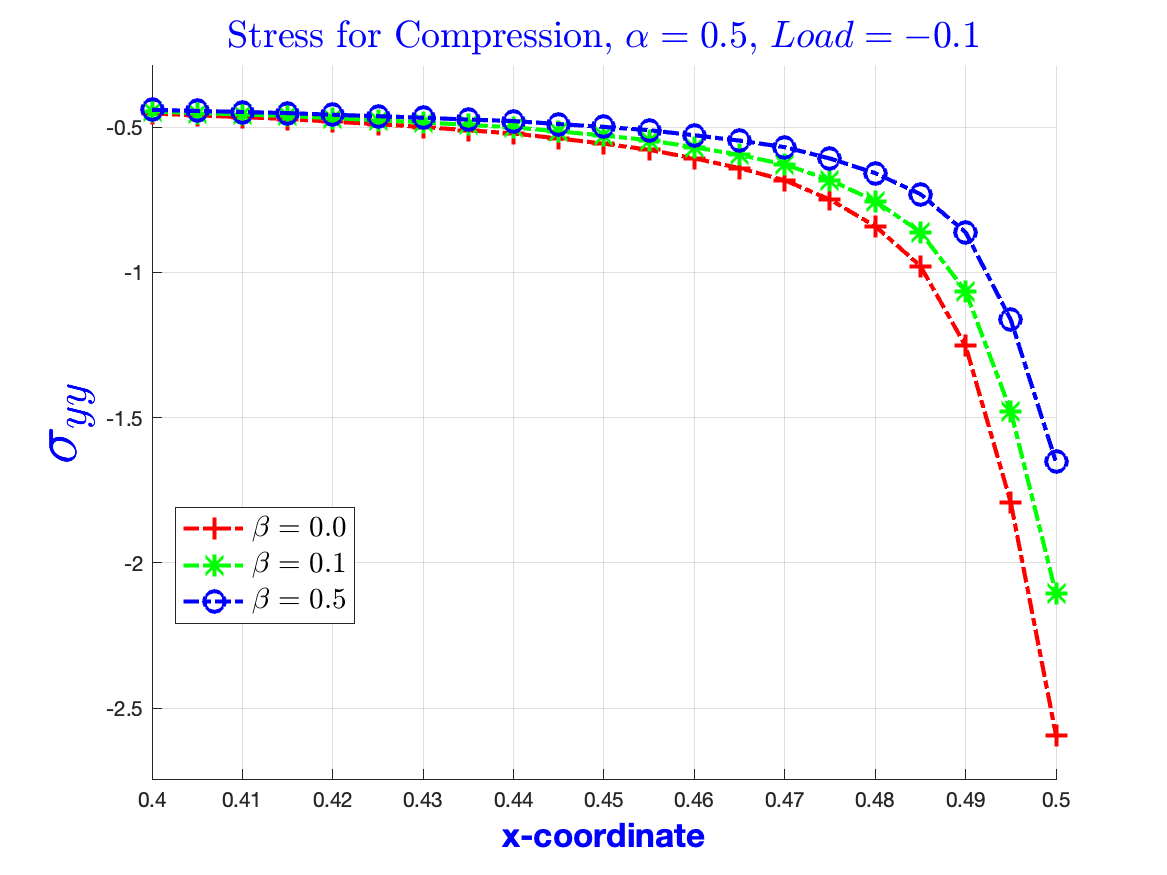}
        \caption{Effect of $\beta$ on $\bfsigma_{yy}$}
        \label{fig:fig2}
    \end{subfigure}
    \vspace{0.5cm} 
    \begin{subfigure}[b]{0.48\textwidth}
        \includegraphics[width=\linewidth]{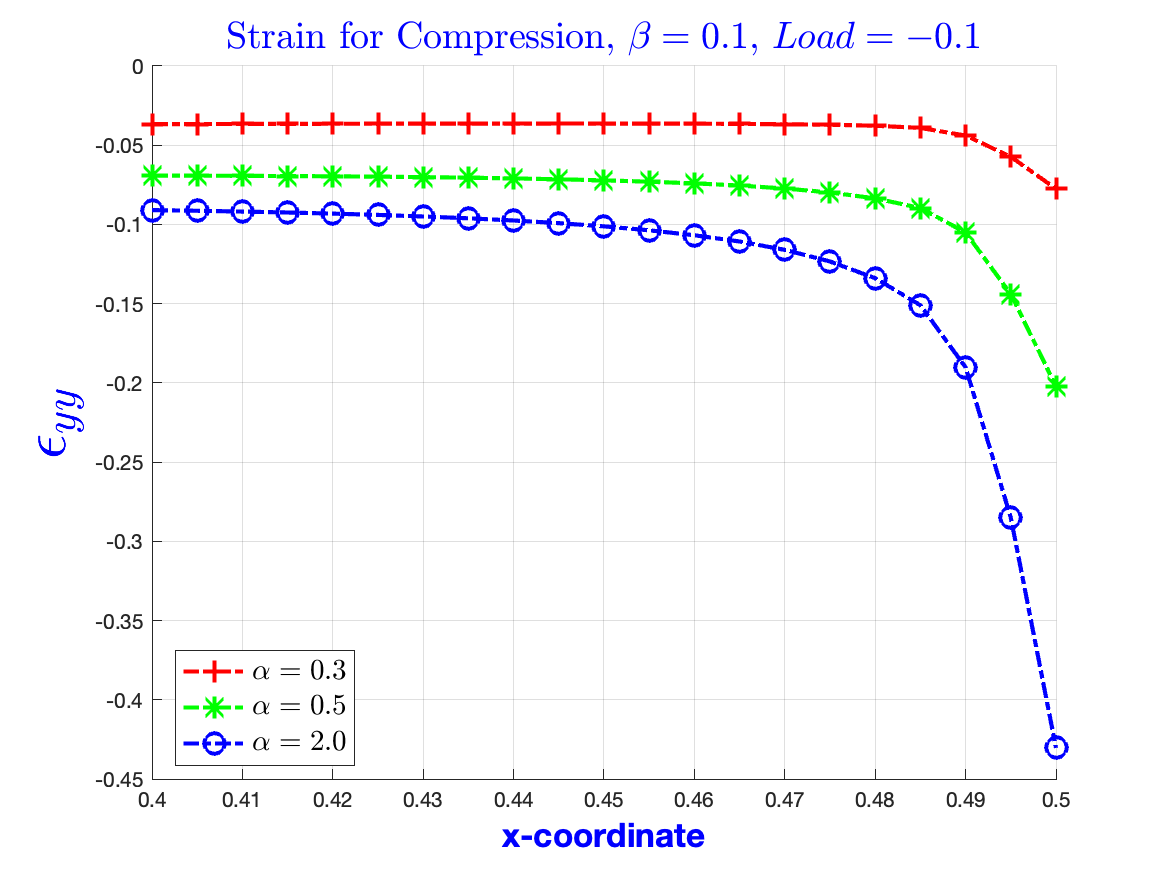}
        \caption{Effect of $\alpha$ on $\bfa{\epsilon}_{yy}$}
        \label{fig:fig3}
    \end{subfigure}
    \hfill 
    \begin{subfigure}[b]{0.48\textwidth}
        \includegraphics[width=\linewidth]{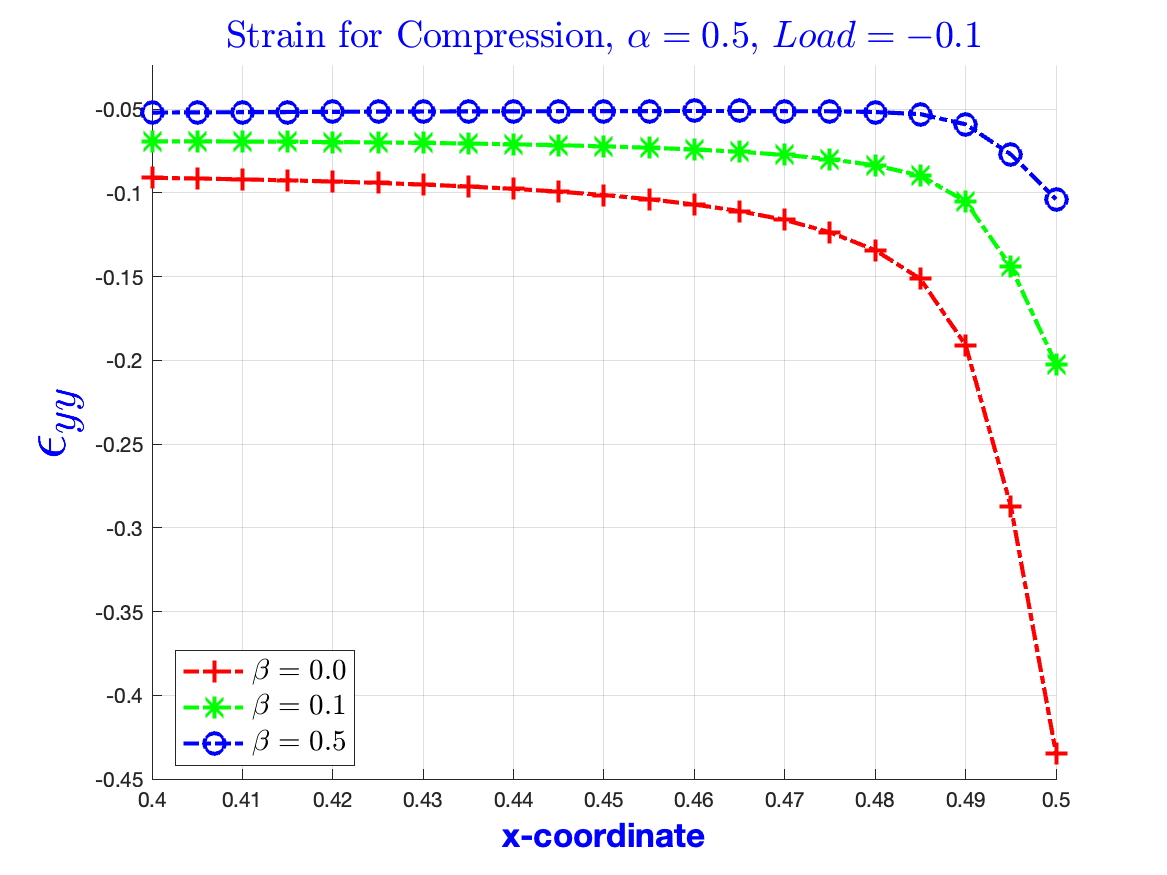}
        \caption{Effect of $\beta$ on $\bfa{\epsilon}_{yy}$}
        \label{fig:fig4}
    \end{subfigure}
    \vspace{0.5cm} 
    \begin{subfigure}[b]{0.48\textwidth}
        \includegraphics[width=\linewidth]{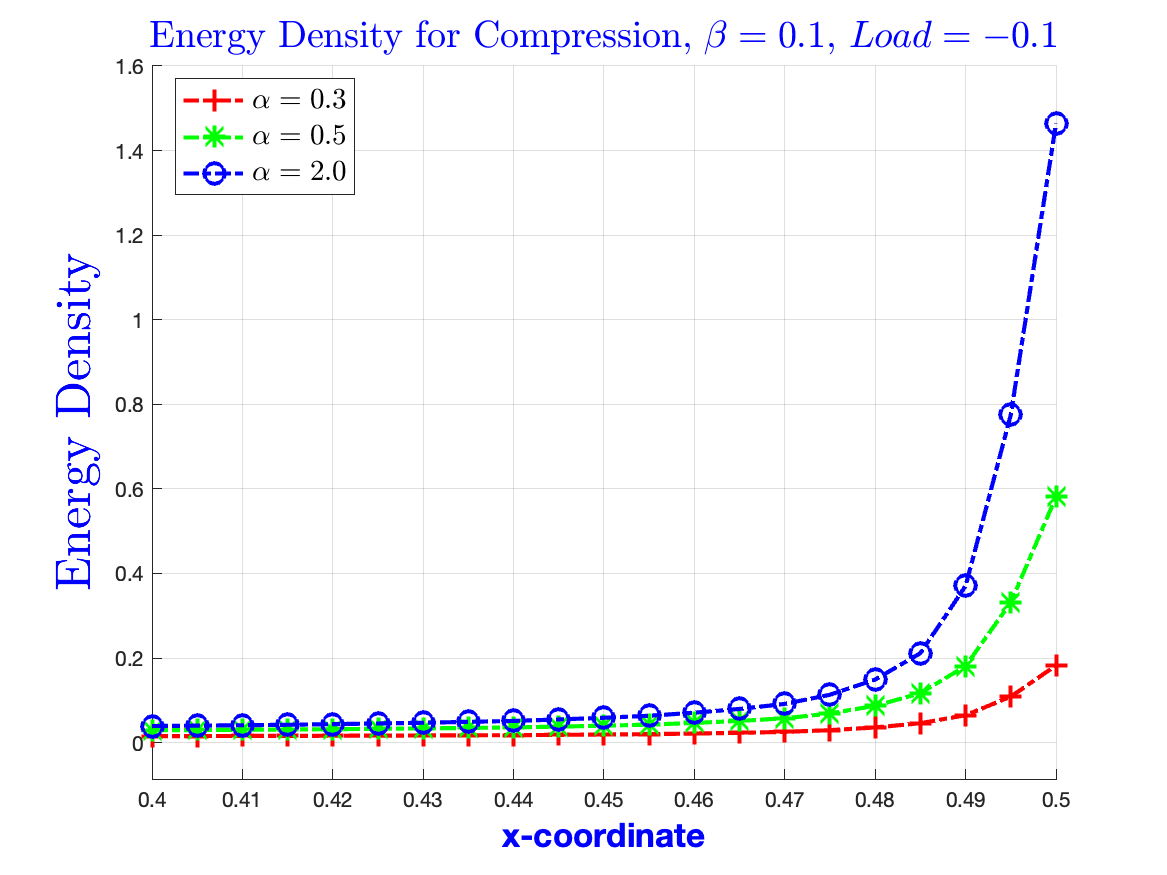}
        \caption{Effect of $\alpha$ on strain-energy density}
        \label{fig:fig5}
    \end{subfigure}
    \hfill 
    \begin{subfigure}[b]{0.48\textwidth}
        \includegraphics[width=\linewidth]{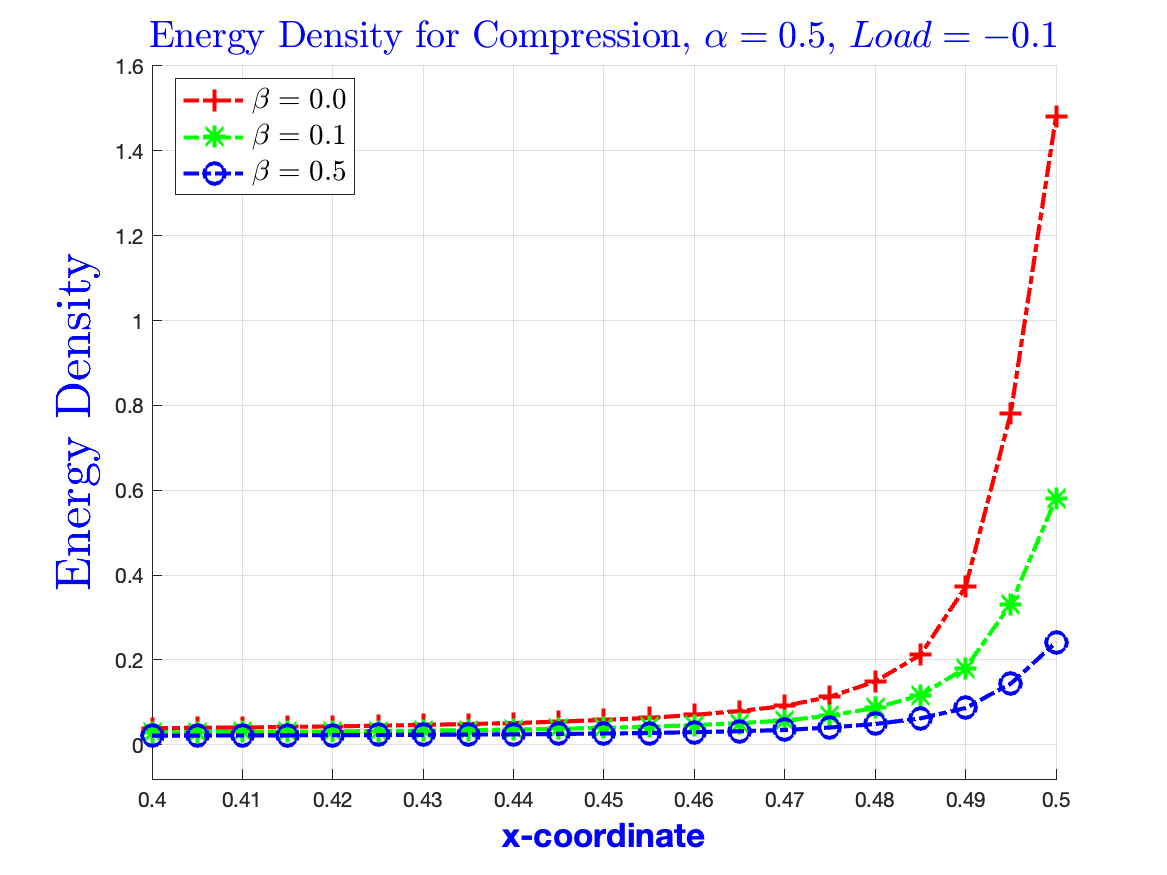}
        \caption{Effect of $\beta$ on strain-energy density}
        \label{fig:fig6}
    \end{subfigure}
\caption{This two-panel figure comprehensively illustrates the impact of the parameters, $\alpha$ and $\beta$, on the crack-tip fields within an orthotropic solid. The material is specifically configured with its fibers aligned orthogonal to the crack.}
    \label{fig:stress_strain_M2}
\end{figure}

Figure~\ref{fig:stress_strain_M2} offers a comprehensive analysis of how the two primary modeling parameters, $\alpha$ and $\beta$, influence the crack-tip fields in a strain-limiting material. The material is modeled as an orthotropic solid with an edge crack, where its internal fibers are aligned perperndicular to the crack. This specific configuration is vital for simulating the behavior of real-world composite materials, like those used in aerospace and automotive industries, where fiber orientation is a key determinant of a component's strength and stiffness. The computations were performed under a fixed top compressive load, allowing for a precise and isolated study of the individual effects of $\alpha$ and $\beta$. The left panel of the figure focuses on the impact of $\alpha$. The results demonstrate that changes in $\alpha$ lead to a proportional and consistent concentration of compressive stress, strain, and strain-energy density around the crack tip. From an engineering perspective, this suggests that $\alpha$ can be used to anchor the model to the principles of classical linear elasticity. By carefully adjusting $\alpha$, engineers can calibrate the model to match the known elastic response of a material, establishing a crucial baseline for more intricate analyses. This capability is essential for validating the model's accuracy and ensuring its predictions are grounded in fundamental material science. Conversely, the right panel reveals the more complex and distinct influence of $\beta$. A key finding here is a paradoxical behavior: as $\beta$ increases, there is a slight decrease of crack-tip stresses,  yet this is coupled with a reduction in both compressive strain and strain-energy density. This is a critical observation for practical engineering applications. It signifies that $\beta$ is the parameter that governs the material's susceptibility to brittle failure under compression. When a material exhibits high stress but low strain and low energy absorption, it has a diminished capacity for plastic deformation. This behavior points to a shift toward a more brittle failure mode, such as sudden crushing or splitting, where intense localized stresses cannot be relieved by material yielding. Thus, $\beta$ provides engineers with a powerful tool for simulating materials where the stress response is highly localized and decoupled from the strain response, enabling more accurate predictions for the behavior of brittle composites or materials under extreme loading conditions. The ability to independently control $\alpha$ and $\beta$ allows for a nuanced and accurate representation of complex material failure mechanisms.

\section{Conclusion}
\label{sec:conclusion}

This study successfully developed and validated a novel mathematical model that employs an algebraically nonlinear constitutive relationship to characterize the mechanical behavior of orthotropic elastic materials. Our primary objective was to provide a deeper understanding of the complex crack-tip fields in strain-limiting bodies, with a specific focus on the influence of distinct fiber orientations. The model's robustness is in the choice of the constitutive response function, which is uniformly buoyant, monotonic, Lipschitz continuous, and coercive, which were fundamental to ensuring the mathematical well-posedness of our continuous Galerkin formulation and, consequently, the existence and uniqueness of the weak solution. By combining the nonlinear constitutive law with the fundamental principle of linear momentum balance, we were able to formulate a comprehensive boundary value problem. To solve this intricate, vector-valued, quasi-linear elliptic problem, we implemented a sophisticated computational approach: a hybrid of Picard's iterative algorithm and the continuous conforming Galerkin finite element method. This enabled us to achieve accurate and efficient numerical solutions, effectively capturing the complex material interactions at the crack tip. Furthermore, our detailed sensitivity analysis revealed the critical impact of key model parameters, varying loading conditions, and different fiber orientations on the material's overall behavior.\\

Our analysis yielded several pivotal findings with significant engineering implications:
\begin{itemize}
    \item \textbf{The Role of $\alpha$ vs. $\beta$}: A key observation from our simulations is the contrasting behavior of the two modeling parameters, $\alpha$ and $\beta$. Increasing the parameter $\alpha$ consistently produced similar effects to a classical linear elastic model. It intensified compressive stress, strain, and strain-energy density at the crack tip, indicating that $\alpha$ can be used to anchor the model to linear elastic theory. In stark contrast, increasing the parameter $\beta$ had a paradoxical and far more complex effect. Higher values of $\beta$ led to a slight decrease in compressive stresses at the crack tip but simultaneously caused a marked reduction in both compressive strain and strain-energy density. This is a critical engineering insight: it suggests that $\beta$ governs the material's susceptibility to brittle failure under compression. When a material can sustain high stress with very little deformation and energy absorption, it has a diminished capacity to tolerate plastic deformation. This behavior points toward a shift from a ductile to a brittle failure mode, where localized high stresses cannot be dissipated, leading to sudden and catastrophic failure, such as crushing or splitting. Therefore, $\beta$ serves as a powerful parameter for engineers to simulate and predict brittle behavior in materials that don't conform to classical elastic-plastic models.

    \item \textbf{Strain Energy Density as a Failure Criterion}: A consistent finding across all our simulations is that the strain energy density is highest in the immediate vicinity of the crack tip, regardless of the values of $\alpha$ and $\beta$. This localized concentration of energy provides a robust local fracture criterion for predicting crack initiation and growth under various loading conditions. From an engineering standpoint, this criterion is invaluable. It allows for the prediction of when a crack will begin to propagate, as the material accumulates energy beyond its critical threshold. This extends the model's applicability from simple mechanical loading to more complex scenarios, such as coupled thermo-mechanical loads. The ability to quantify this energy distribution offers crucial insights into a material's ductility or brittleness at the point of failure, which is essential for developing accurate predictive models for component integrity and lifespan.
\end{itemize}

This study establishes a strong foundation that can be significantly expanded to further advance our understanding of crack behavior in complex materials.
\begin{enumerate}
    \item \textbf{Porous and cellular materials}: A critical next step is to investigate crack-tip fields within porous elastic bodies where material properties are intrinsically linked to density. This would enable more accurate modeling of lightweight materials with internal voids, which are increasingly common in modern engineering \cite{yoon2024finite,gou2025computational}.
    \item \textbf{Three-dimensional analysis}: Moving beyond two-dimensional simplifications is essential for a more realistic representation of crack propagation and stress distribution in actual components. Resolving crack-tip fields in three-dimensional bodies will provide insights vital for designing safer and more durable structures \cite{mallikarjunaiah2025crack,gou2023computational,gou2023finite}.
    \item \textbf{Computational methodologies}: A thorough numerical analysis comparing the performance of both continuous and discontinuous Galerkin finite element methods for these types of problems is another important research topic. Such an analysis would optimize computational efficiency, enhance the accuracy of numerical predictions, and contribute significantly to the field of computational fracture mechanics \cite{manohar2024hp}.
   \item \textbf{Quasistatic crack propagation under compressive loading}: A logical and significant extension of our current model is to study quasistatic crack propagation under compressive loading \cite{zou2021combined}. This is a crucial area of research, as it moves beyond analyzing a static crack to simulating how a crack actually grows and evolves over time. We can formulate this problem within our existing framework by using the Francfort-Marigo energy minimization principle \cite{francfort1998revisiting,manohar2025convergence,manohar2025adaptive,lee2022finite,yoon2021quasi,fernando2025xi,fernando2025textsf}. This principle describes a material's total energy, ${E}$, as a function of both the displacement field and the crack geometry. 
\begin{equation}\label{eq:tenergy}
{E}(\bfa{u}, \; \Gamma_c):=\int_{\Lambda \setminus \Lambda_c} \mathcal{W}(\bfa{u})  \; d\bfx + \mathcal{G}_c \mathcal{H}^{1}(\Gamma_c),
\end{equation}
where $\mathcal{W}(\cdot) \colon H^{1}(\Lambda) \to \mathbb{R}$ is the elastic energy, $\bfa{u} \colon \Lambda \to \mathbb{R}^2$ is the displacement field, $\mathcal{H}^{1}$ denotes Hausdorff measure, and $\mathcal{G}_c$ denotes the \textit{critical energy release rate} (or \textit{fracture toughness}) of the material. The total energy consists of two parts: the elastic energy stored in the material, $\mathcal{W}(\bfa{u})$, and the energy required to create a new crack surface, $\mathcal{G}_c \mathcal{H}^{1}(\Lambda_c)$. By minimizing this total energy, our model can predict the direction and conditions under which a crack will grow, providing a powerful tool for understanding and predicting material failure under compressive loads.

\end{enumerate}

\section{Acknoledgement}
SG would like to thank the University of Texas Rio Grande Valley for providing a Presidential Research Fellowship during his PhD studies. SMM's work is supported by the National Science Foundation under Grant No. 2316905.
 
\bibliographystyle{plain}
\bibliography{ref}

\end{document}